\documentclass[10pt]{amsart}                          
\usepackage{epsfig}
\usepackage{amssymb,amscd}     
\usepackage{amsthm,amsgen,amsmath,epic,psfrag}

\newcommand{\R}{{\mathbb R}}  
\newcommand{\I}{{\mathbb I}}                        
\newcommand{\C}{C_n({\mathbb R}^m)}

\newcommand{\la}{\langle [}                          
\newcommand{\ra}{] \rangle}                          
                          
\newcommand{\und}{\underline} 
\newcommand{\n}{\underline{n}}

\newcommand{\wt}{\widetilde}                          
\newcommand{\comment}[1]{}%{{\bf #1}}%{} 
\newcommand{\skiphome}[1]{#1}                         
                         
\theoremstyle{plain}                          
\newtheorem{theorem}{Theorem}[section]                          
\newtheorem{proposition}[theorem]{Proposition}                          
\newtheorem{lemma}[theorem]{Lemma}                          
\newtheorem{corollary}[theorem]{Corollary}                          
 
\newtheorem{mydiagram}[theorem]{Figure}

\theoremstyle{definition}                          
\newtheorem{definition}[theorem]{Definition}                          
\newtheorem*{example}{Example}                          
                         
\theoremstyle{remark}                          
\newtheorem*{remark}{Remark}                          
\newtheorem*{notation}{Notation}
\newtheorem*{convention}{Convention}                                                    

\newcommand{\refT}[1]{Theorem~\ref{T:#1}}                          
\newcommand{\refC}[1]{Corollary~\ref{C:#1}}                          
\newcommand{\refP}[1]{Proposition~\ref{P:#1}}                          
\newcommand{\refD}[1]{Definition~\ref{D:#1}}                          
\newcommand{\refL}[1]{Lemma~\ref{L:#1}}

\begin{document}                                                  

\title{Manifold-theoretic compactifications of configuration spaces}       
\author{Dev P. Sinha}       
\address{Department of Mathematics, University of Oregon, Eugene, OR
97403}       
\email{dps@math.uoregon.edu}       

       
\begin{abstract}       
We present new definitions for and give a comprehensive treatment  of the
canonical compactification of configuration spaces due to  Fulton-MacPherson
and  Axelrod-Singer in the setting of smooth manifolds, as well as a
simplicial variant of this compactification initiated by Kontsevich.  Our constructions are
elementary and give simple global coordinates for the compactified
configuration space of a general manifold embedded in Euclidean space.   We
stratify the canonical compactification, identifying the 
diffeomorphism types of the strata in terms of spaces of  configurations in
the tangent bundle, and give completely explicit local coordinates  around
the strata as needed to define a manifold with corners.    We analyze the quotient map from the canonical to the
simplicial  compactification, showing it is a homotopy equivalence.  Using global
coordinates we define projection maps and diagonal maps, which for the simplicial variant
satisfy cosimplicial identities.  
\end{abstract}

\maketitle                     

{\tableofcontents}

\section{Introduction} 

Configuration spaces are fundamental objects of study in geometry  and
topology, and over the past ten years, functorial compactifications of
configuration  spaces have been an important technical tool.  We review the
state of this active  area after giving our definitions.

\subsection{Basic definitions}

We first choose compact notation to manage products of spaces.

\begin{notation}
If $S$ is a finite set, $X^S$ is the product $X^{\#S}$ where $\# S$
is the cardinality of $S$.  Consistent with this, if $\{X_s\}$ is a
collection of spaces indexed by $S$, we let $(X_s)^S = \prod_{s \in S}
X_s$.  For coordinates in either case we use
$(x_s)_{s\in S}$ or just $(x_s)$ when $S$ is understood. Similarly,  a
product of maps $\prod_{s\in S} f_s$ may be written$(f_s)_{s \in S}$ or
just $(f_s)$. We let $\n$ denote the set $\{1, \ldots, n\}$, our most
common indexing set.
\end{notation}

\begin{definition}
If $M$ is a smooth manifold, let $C_n(M)$ be the subspace of  
$(x_i) \in M^{\n}$ such that $x_i \neq x_j$ if $i \neq j$.  Let
$\iota$ denote the inclusion of $C_n(M)$ in $M^{\n}$.   
\end{definition}

Suppose that $M$ were equipped with a metric.  The main
compactification  which we study, $C_n[M]$, is homeomorphic to the
subspace of $C_n(M)$ for  which $d(x_i, x_j) \geq \epsilon$ for some
sufficiently small $\epsilon$.   From this model, however, it is not clear
how $C_n(M)$ should be a subspace of the compactification, much less how to
establish functorality or more delicate properties that we will develop.  

\begin{definition} 
For $(i,j) \in C_2(\n)$,
let $\pi_{ij} \colon \C \to S^{m-1}$ be the map which sends $(x_i)$
to the unit vector in the direction of $x_i - x_j$. 
Let $I$ be the closed interval from $0$ to $\infty$, the one-point
compactification of $[0, \infty)$.  For
$(i,j,k) \in C_3(\und{n})$ let $s_{ijk} \colon \C \to I = [0, \infty]$ 
be the map which sends $(x_i)$ to 
$\left(|x_i - x_j|/|x_i - x_k|\right)$. 
\end{definition} 

Our compactifications are defined as closures, for which we also 
set notation.

\begin{notation}
If $A$ is a subspace of $X$, we let $cl_X(A)$, or simply $cl(A)$ if by
context $X$ is understood, denote the closure of
$A$ in $X$.  
\end{notation}

From now on by a manifold $M$ we mean a submanifold of some $\R^m$, so that 
$C_n(M)$ is a submanifold of $C_n(\R^m)$.  
For $M = \R^m$, we specify that $\R^m$ is a submanifold of itself through 
the identity map.

\begin{definition}\label{D:maindef}
Let $A_n[M]$, the main ambient space in which we work, be the product  
$M^{\und{n}} \times  (S^{m-1})^{C_2(\n)}
\times I^{C_3(\n)},$ and similarly
let $A_n\la M \ra = M^{\n} \times (S^{m-1})^{C_2(\n)}$.
Let $$\alpha_n = \iota \times \left(\pi_{ij}|_{C_n(M)}\right) \times
\left((s_{ijk})|_{C_n(M)}\right) : C_n(M) \to A_n[M]$$ 
and define $C_n[M]$ to be $cl_{A_n[M]}\left(im(\alpha_n)\right)$.  Similarly, let
$\beta_n = \iota \times \left(\pi_{ij}|_{C_n(M)}\right) :
 C_n(M) \to A_n\la M \ra$ and define $C_n \la M \ra$ to be 
$cl_{A_n\la M \ra}\left(im(\beta_n)\right)$.
\end{definition}

We will show in \refT{mfld} that 
$C_n[M]$ is a manifold with corners whose diffeomorphism type
depends only on that of $M$.  Because
$A_n[M]$ is compact when $M$ is and $C_n[M]$ is closed in
$A_n[M]$, we immediately have the following.  

\begin{proposition} 
If $M$ is compact, $C_n[M]$ is compact. 
\end{proposition} 

We call $C_n[M]$ the canonical compactification of $C_n(M)$ and $C_n \la M
\ra$ the simplicial variant.  When $M$ is not compact but is equipped with a complete metric, it is natural to call $C_n[M]$ the canonical completion of $C_n(M)$.

\subsection{Review of previous work}

The compactification $C_n[M]$ first appeared in work of Axelrod and
Singer \cite{Axel94}, who translated the
definition of Fulton and MacPherson in  \cite{Fult94} as a closure in a
product of blow-ups from algebraic geometry to the setting of manifolds
using spherical blow-ups.  
Kontsevich made  similar constructions at about the same time as
Fulton and MacPherson, and his later definition in 
\cite{Kont99} coincides with our $\wt{C}_n \la \R^m\ra$,
although it seems that he was trying to define $\wt{C}_n[\R^m]$. 
Kontsevich's oversight was corrected in \cite{Gaif03}, in which Gaiffi
gives a definition of $C_n[\R^m]$ similar to ours, generalizes the
construction for arbitrary hyperplane arrangments over the real numbers,  
gives a pleasant
description of the category of strata using the language of blow-ups of
posets from \cite{FeKo03}, and also treats blow-ups for stratified spaces
locally and so gives rise to a new definition of
$C_n[M]$, but one which is less explicit than ours and thus less suited for
the applications we develop.  An alternate
approach to $C_n[M]$ through the theory of operads as pioneered 
by Getzler and Jones \cite{GeJo95} was fully developed and
extended to arbitrary manifolds
by Markl \cite{Mark99}.
 
Axelrod and Singer used these compactifications to define 
invariants of three-manifolds coming
from Chern-Simons theory, and these constructions have generally
been vital in quantum
topology \cite{Bott94, KuTh01, BGRT02, Poir02}.   Extensive use of similar 
constructions has been made in the setting of hyperplane
arrangments \cite{DCPr95, Yuzv97} over the complex numbers.  These
compactifications have also inspired new computational 
results \cite{Kriz95, Tota96}, and they canonically realize the homology
of $C_n(\R^m)$ \cite{CMS03}.  We came to the present definitions of these
compactifications so we could define
maps and boundary conditions needed for applications to knot
theory \cite{BCSS03, Sinh02}.

New results in this paper include full proofs of many folk theorems, and the following:

\begin{itemize}
\item A construction for general manifolds which bypasses the need for
blow-ups, uses simple global coordinates, and through which
functorality is immediate.
\item Explicit description of the strata in terms of spaces of 
configurations in the tangent bundle.
\item Coordinates about strata which may easily be used for 
transversality arguments.
\item Full treatment of the simplicial variant, including a proof that
the projection from the canonical compactification to the simplicial
one is a homotopy equivalence.
\item A clarification of the central role which Stasheff's associahedron
plays in this setting.
\item Explicit identification of these compactifications as subspaces
of familiar spaces.
\item Constructions of diagonal maps, projections, and substitution maps 
as needed for applications.  \end{itemize}

In future work \cite{McSS03}, we will use these constructions to define an operad
structure on these compactifications of configurations in Euclidean space, which
has consequences in knot theory.  This operad structure was first applied in
\cite{GeJo95}.

We also hope that a unified and explicit exposition of these
compactifications using our simplified definition  could be of help, especially to
those who are  new to the subject.

\subsection{A comment on notation, and a little lemma}
 
There are two lines of notation for configuration spaces of manifolds
in  the literature, namely $C_n(M)$ and $F(M,n)$.  Persuaded by Bott,
we  choose to use the $C_n(M)$ notation.  Note, however, that $C_n(M)$
in  this paper is $C^0_n(M)$ in \cite{Bott94} and that $C_n[M]$ in this
paper is $C_n(M)$ in \cite{Bott94}.   Indeed, we warn the reader to pay
close  attention to the parentheses in our notation: 
$C_n(M)$ is the open configuration space;  $C_n[M]$ is the
Fulton-MacPherson/Axelrod-Singer compactification, its canonical completion;
$C_n\la M \ra$, the simplicial variant, is a quotient of $C_n[M]$;
$C_n\{M\}$, an auxilliary construction, is a subspace of $C_n[M]$
containing only one additional stratum.  We suggest that those who  choose
to use $F(M,n)$  for the open configuration space use $F[M,n]$ for the
compactification. 

As closures are a central part of our definitions, we need a lemma from
point-set topology that open maps commute with taking closures.

\begin{lemma}\label{L:projclose}
Let $A$ be a subspace of $X$, and let $\pi : X \to Y$ be an 
open map.  Then $\pi(cl_X(A)) \subseteq cl_Y(\pi(A))$.  If $cl_X(A)$ is
compact (for example, when $X$ is) then this inclusion is an equality.
\end{lemma}

\begin{proof}
First, $\pi^{-1}(cl_Y(\pi(A)))$ is closed in $X$ and
contains $A$, so it contains $cl_X(A)$ as well.  Applying
$\pi$ to this containment we see that $\pi(cl_X(A)) \subseteq cl_Y(\pi(A))$.

If $cl_X(A)$ is compact, so is $\pi(cl_X(A))$, which is thus closed in $Y$.
It contains $\pi(A)$, therefore $cl_Y(\pi(A)) \subseteq \pi(cl_X(A))$.
\end{proof}

\subsection{Acknowledgements} 
The author would like to thank Dan Dugger for providing a proof and references 
for \refL{push},  Tom Goodwillie for providing the main idea for  \refL{3d4c},  
Ismar Volic for working with the author
on  an early draft of this paper, Matt Miller for a careful reading, and 
Giovanni Gaiffi and Eva-Maria Feitchner for sharing preprints of their work.

\section{A category of trees and related categories} 

In order to understand the compactifications $C_n[M]$ we have to 
understand their strata, which are naturally labelled by a poset (or 
category) of trees.  

\begin{definition}  
Define an $f$-tree to be a rooted, connected tree, with labelled
leaves, and with no bivalent  internal vertices.  Thus, an $f$-tree $T$ 
is a connected acyclic graph with a specified vertex $v_0$ called the root. 
The root may have any valence, but other vertices may not be bivalent. 
The univalent vertices other than perhaps the root are called leaves, and 
each leaf is labelled uniquely with an element of $\und{\#l(T)}$, where 
$l(T)$ is the set of leaves of $T$ and $\#l(T)$ is its cardinality.  
\end{definition}

\begin{center}
\begin{minipage}{10cm}
\begin{mydiagram}\label{F:T}\begin{center}
A tree $T$.
{
\psfrag{1}{1}  \psfrag{2}{2}  \psfrag{3}{3}  \psfrag{4}{4}  \psfrag{5}{5}  \psfrag{6}{6}  \psfrag{7}{7}
\psfrag{A}{$v_1$}
\psfrag{B}{$v_0$}
\psfrag{C}{$v_2$}
\psfrag{D}{$v_{3}$}
$$\includegraphics[width=3.5cm]{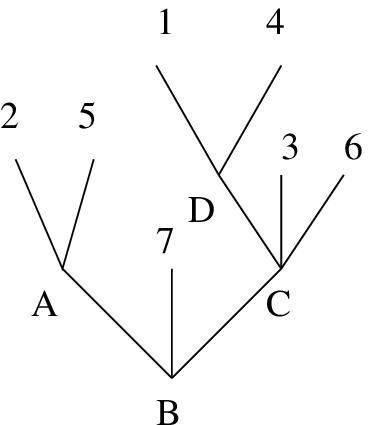}$$
}

\end{center}
\end{mydiagram}
\end{minipage}
\end{center}

In an $f$-tree there is a unique path from any vertex or edge to the root 
vertex, which we call its root path.  We say that one vertex or edge lies 
over another if the latter is in the root path of the former.
For any edge, its boundary vertex closer to the root is called its
initial vertex, and its other vertex is called its terminal vertex.  If two
edges share the same initial vertex, we call them coincident.  For
a vertex $v$ there is a canonical ordering of edges for which $v$ is
initial, the collection of which we call $E(v)$, the group of edges
coincident at $v$.  Namely, set $e < f$ if the smallest label for a leaf over
$e$ is smaller than that over $f$.  We may use this ordering to name these
edges $e_1(v), \ldots, e_{\#v}(v)$, where $\#v$ is the number of edges 
in $E(v)$. 

We will be interested in the set of $f$-trees as a set of objects in a 
category in which morphisms are defined by contracting edges. 

\begin{definition} 
Given an $f$-tree $T$ and a set of non-leaf edges $E$ the 
contraction of $T$ by $E$ is the tree $T'$ obtained by
taking each edge $e \in E$, identifying its initial vertex with its
terminal vertex, and then removing $e$ from the set of edges. 
\end{definition} 

\begin{definition}  
Define $\Psi_{\und{n}}$ to be the category whose objects are $f$-trees with
$n$  leaves.  There is a (unique) morphism in $\Psi_{\und{n}}$ from $T$ to
$T'$ if $T'$ is  isomorphic to a contraction of $T$ along some set of edges. 
\end{definition} 

\skiphome{
\begin{center}
\begin{minipage}{10cm}
\begin{mydiagram}\label{F:Psi3}\begin{center}
The category $\Psi_3$.
{
$$\includegraphics[width=6cm]{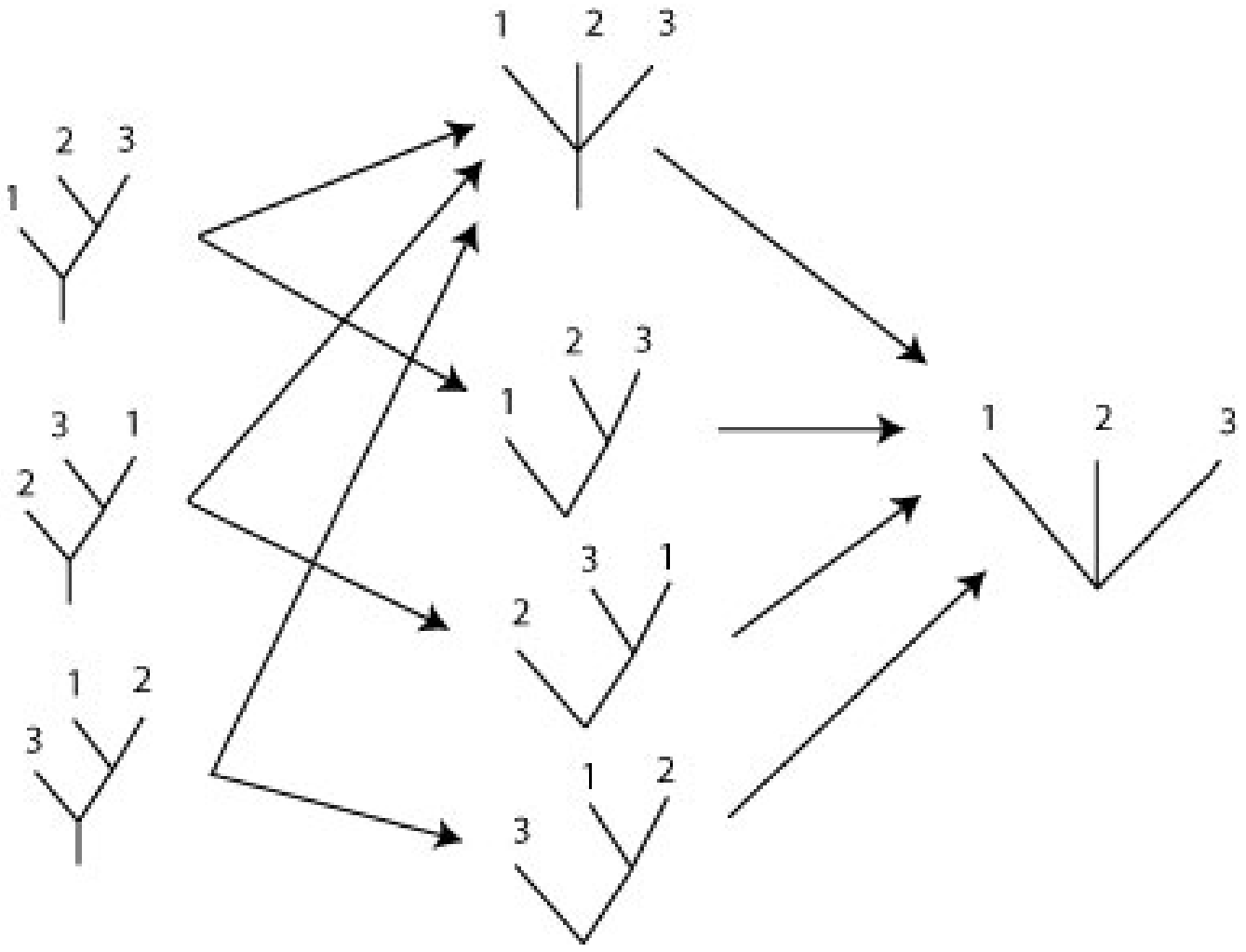}$$
}

\end{center}
\end{mydiagram}
\end{minipage}
\end{center}

%\begin{figure}
%$$\includegraphics[width=8cm]{Psi3.eps} $$
%\caption{The category $\Psi_3$.}
%\end{figure}
}

Finally, let $V(T)$ denote the set of non-leaf vertices of $T$. 
Let $V^i(T)$ denote its subset of internal vertices (thus only 
excluding the root).  Note that a morphism in
$\Psi_{\und{n}}$ decreases the number of internal vertices, which is zero for
the terminal object in $\Psi_{\und{n}}$.  Let $\wt{\Psi_n}$ be the full subcategory 
of $f$-trees whose root is univalent (informally, trees with a trunk).  Note that
$\wt{\Psi_n}$ has an operad structure, as defined in \cite{GeJo95, MSS02} 
(but $\Psi_n$ does
not have one since we do not allow bivalent vertices).

\medskip

It is useful to have facility with categories that are essentially 
equivalent to  $\Psi_{\und{n}}$.  We will define these categories through
the notions of parenthesization and exclusion relation.  Further
equivalent constructions include the collections 
of screens of Fulton and MacPherson \cite{Fult94}.  
The best perspective on these categories 
is given by the combinatorial blow-up of
Feitchner and Kozlov \cite{FeKo03}.  Indeed, Gaiffi shows in \cite{Gaif03}
that  the poset of strata of a blow-up of an arrangment is
the combinatorial blow-up of the orignial poset associated to the
arrangment.  Since we focus not on general blow-ups but on compactified configuration
spaces in particular, we choose more concrete manifestations of this category.

\begin{definition} 
A (partial)  parenthesization $\mathcal{P}$ of a set $S$ is a 
collection $\{ A_\alpha \}$ of nested subsets of $S$, 
each of cardinality greater than one. By nested we mean that,
for any $\alpha,\beta$, the intersection $A_\alpha \cap A_\beta$ is 
either $A_\alpha$, $A_\beta$ or empty.  The 
parenthesizations of $S$ form a poset, which we call $Pa(S)$, in which 
$\mathcal{P} \geq \mathcal{P}'$ if $\mathcal{P} \subseteq \mathcal{P}'$. 
\end{definition} 

Parenthesizations are related to trees in that they may keep track
of sets of leaves which lie over the vertices of a tree. 

\begin{definition} 
Define $f_1 \colon \Psi_{\und{n}} \to Pa(\und{n})$ by sending a tree $T$ 
to the  collection of sets $\{A_v\}$, where $v \in V^i(T)$ and $A_v$ is the
set of indices of leaves which lie over $v$. 
Define $g_1 \colon Pa(\und{n}) \to \Psi_{\und{n}}$ by sending a
parenthesization to a tree with the following data
\begin{enumerate} 
\item One internal vertex $v_\alpha$ for each $A_\alpha$.
\item An edge between $v_\alpha$ and $v_\beta$ if 
$A_\alpha \subset A_\beta$ but there is no proper 
$A_\alpha \subset A_\gamma \subset A_\beta$.
\item A root vertex with an edges connecting it to each internal vertex 
corresponding to a maximal $A_\alpha$.
\item Leaves with labels in $\und{n}$ with an edge connecting the $i$th 
leaf  to either the vertex $v_\alpha$ where $A_\alpha$ is the minimal 
set  containing $i$, or the root vertex if there is no such $A_\alpha$.
\end{enumerate}
\end{definition} 

We leave to the reader the straightforward verification that $f_1$ and 
$g_1$ are well defined and that the following proposition holds. 

\begin{proposition} 
The functors $f_1$ and $g_1$ are isomorphisms between the categories 
$\Psi_{\und{n}}$ and $Pa(\und{n})$. 
\end{proposition} 

Another way in which to account for the data of which leaves lie above 
common vertices in a tree is through the notion of an exclusion
relation. 

\begin{definition} 
Define an exclusion relation $R$ on a set $S$ to be a subset of $C_3(S)$ 
such that the following properties hold.
\begin{enumerate} 
\item If $(x,y),z \in R$, then $(y,x),z \in R$ and $(x,z),y \notin R$.
\item If  $(x,y), z \in R$ and $(w,x), y \in R$, then $(w,x),z \in R$. 
\end{enumerate} 
If $(x,y),z \in R$ we say that $x$ and $y$ exclude $z$.
Let $Ex(S)$ denote the poset of exclusion relations on $S$, where the 
ordering is defined by inclusion as subsets of $C_3(S)$. 
\end{definition} 

We now construct exclusion relations from parenthesizations, and 
vice-versa. 

\begin{definition} \label{D:extree}
Let $f_2 \colon Pa(\und{n}) \to Ex(\und{n})$ be defined by setting 
$(i, j),k \in R$ if $i,j \in A_\alpha$ but $k \notin A_\alpha$ for some 
$A_\alpha$ in the  given parenthesization.  Define $g_2 \colon
Ex(\und{n}) \to Pa(\und{n})$  by, given an exclusion relation $R$, taking
the collection of sets $A_{\sim i,\neg k}$ where $A_{\sim i, \neg k}$ is
the set of all $j$ such that $(i,j), k \in R$, along with $i$ when there is
such a $j$.  Let $Tr = g_1 \circ g_2 : Ex(\und{n}) \to \Psi_{\und{n}}$ and
let ${\mathcal{E}x} = f_2 \circ f_1$.
\end{definition} 

As above, we leave the proof of the following elementary proposition to the 
reader.

\begin{proposition}
The composite $f_2 \circ g_2$ is the identity functor.  If $f_2(\mathcal{P}) 
= f_2(\mathcal{P'})$, then $\mathcal{P}$ and $\mathcal{P'}$ may only
differ by whether or not they contain the set $\und{n}$ itself.
\end{proposition} 

\section{The stratification of the basic compactification}\label{S:strata}

This section is the keystone of the paper. 
We first define a stratification of $C_n[M]$ through coordinates 
as a subspace of $A_n[M]$. 
For our purposes, a stratification is any expression of a 
space as a finite disjoint union of locally closed
subspaces called strata, which are usually manifolds, such that the
closure of each stratum is its union with other strata.  We 
will show that when $M$ has no boundary, the stratification we
define through coordinates coincides with
the stratification of $C_n[M]$ as a manifold with corners.   
The strata of $C_n[M]$ are individually simple to describe, so constructions and maps on 
$C_n[M]$ are often best understood in terms of these strata.   

Before treating $C_n[M]$ in general, we would like to be completely
explicit about the simplest possible case, essentially $C_2[\R^m]$.

\begin{example}
Let $C_2^*(\R^m) \cong \R^m- 0$ be the subspace of points $(0, x \neq 0) 
\in C_2(\R^m)$ and consider it as the subspace of $\R^m\times
S^{m-1}$ of points $(x \neq 0,  \frac{x}{||x||})$.  The
projection of this subspace onto $S^{m-1}$ coincides with the tautological
positive ray bundle over $S^{m-1}$, which is a trivial bundle.  The closure
$C_2^*[\R^m]$ is the non-negative ray bundle, which is diffeomorphic to
$S^{m-1} \times [0, \infty)$.  Projecting this closure onto $\R^m$ is a
homeomorphism when restricted to
$\R^m- 0$, and the preimage of $0$ is a copy of $S^{m-1}$, the stratum
of added points.  
Thus, $C^*_2[\R^m]$ is diffeomorphic to 
the blow-up of $\R^m$ at $0$, in which one replaces $0$
by the sphere of directions from which it can be approached.   
Through this construction, $C^*_2[\R^m]$ has simple global 
coordinates inherited from $\R^m\times S^{m-1}$.
\end{example}

\subsection{Stratification of $C_n[M]$ using coordinates in $A_n[M]$}

We proceed to define a stratification for $C_n[M]$ by associating an 
$f$-tree to each point in $C_n[M]$.
 
\begin{definition}\label{D:Tx}
Let $x =  \left(( x_i), ( u_{ij}), (d_{ijk}) \right) \in C_n[M]$.  Let 
$R(x)$ be the exclusion relation defined by $(i,j), k \in R(x)$ if 
$d_{ijk} = 0$.  
Let $T(x)$ be equal to either $Tr(R(x))$ or, if all of the $x_i$ are 
equal, the $f$-tree obtained by adding a new root to $Tr(R(x))$.
\end{definition}  

Note that because $d_{ijk} d_{i \ell j} = d_{i \ell k}$ for points in
the image of $C_n(M)$, by continuity this is true for all of $C_n[M]$.  So
if $d_{ijk} = 0 = d_{i \ell j}$, then $d_{i \ell k} = 0$.  Therefore,
$R(x)$ satisfies the last axiom for an  exclusion relation.  The other
axiom is similarly straightforward to check to see that $R(x)$ is
well defined.  

\begin{definition} 
Let $C_T(M)$ denote the subspace of all $x \in C_n [M]$ such that $T(x)
= T$, and let $C_T[M]$ be its closure in $C_n[M]$. 
\end{definition} 

The  following proposition, which gives a first indication of how the 
$C_T(M)$ fit together, is an immediate consequence of the definitions above. 

\begin{proposition} 
Let $s = \{(x_i)_j \}_{j=1}^\infty$ be a
sequence  of of points in $C_n(M)$ which converges to a point in $C_n[M] 
\subset A_n[M]$.   The limit of $s$ is in $C_T[M]$ 
if and only if the  limit of $d(x_i, x_j)/ d(x_i, x_k)$ approaches zero 
for every $(i,j),k \in Ex(T)$ and, in the case where the root valence 
of $T$ is one, we also have that all of the $x_i$ approach the same 
point in $M$. 
\end{proposition} 

To a stratification of a space, one may associate a poset in which 
stratum $\alpha$ is less than stratum $\beta$ if $\alpha$ is contained
in the closure of $\beta$.  

\begin{theorem}\label{T:strata1} 
The poset associated to the stratification of $C_n[M]$ by the $C_T(M)$
is isomorphic to $\Psi_{\und{n}}$. 
\end{theorem} 

\begin{proof} 
This theorem follows from the preceding proposition and the
fact that  if $T \to T'$ is a morphism in $\Psi_{\und{n}}$, then $R(T')$ is
contained in $R(T)$. 
\end{proof}  

\subsection{Statement of the main theorem}

Having established an intrinsic definition for the $C_T(M)$ and a 
combinatorial description of how they fit together, we now set ourselves
to the more difficult task of identifying these spaces explicitly.
We describe the spaces $C_T(M)$ in terms of ``infinitesimal
configurations''.   We will use the term scaling to refer to the action of positive
real numbers on a vector space through scalar multiplication.

\begin{definition}\label{D:ic} 
\begin{enumerate}
\item  Let $Sim_k$ be the subgroup of the group of affine transformations in $\R^m$ 
generated by translation and scalar multiplication.
\item Define $IC_i(M)$ to be the space of $i$ distinct points in $TM$ all
lying in one fiber, modulo the action of $Sim_k$
in that fiber.    Let $p$ be the projection of $IC_i(M)$ onto $M$. 
\end{enumerate}
\end{definition} 

For example $IC_2(M)$ is diffeomorphic to $STM$, the unit tangent
bundle of $M$.   We sometimes refer to $IC_i(M)$ as the space of infinitesimal 
configurations of $i$ points in $M$.

Let $e \in E_0 = E(v_0)$ be a root edge of an $f$-tree
$T$, and let $V(e) \subseteq V^i(T)$ be the set of
internal vertices which lie over $e$. 

\begin{definition} \label{D:DT}
\begin{enumerate}
\item Define $IC_e(M)$ to be subspace of the product
$(IC_{\#v}(M))^{V(e)}$ of tuples of infinitesimal configurations all sitting over
the same point in $M$.  
\item Let $p_e$ be the map from
$IC_e(M)$ onto $M$ defined projecting onto that point. 
\item Let $D_T(M)$ be the subspace of $(IC_e(M))^{E_0}$ of points 
whose image under $(p_e)$ in $(M)^{E_0}$ sits in $C_{\#v_0}(M)$.
\end{enumerate}
\end{definition}

In other words, a point in $D_T(M)$ is a collection of $\#v_0$ distinct
points $(x_e)_{e \in E_0}$ in $M$ with a collection of $\#v(e)$ 
infinitesimal configurations at each $x_e$.

\skiphome{
\begin{center}
\begin{minipage}{10cm}
\begin{mydiagram}\label{F:pictDT}\begin{center}
A  point in $D_T(M)$ with $T$ from Figure~\ref{F:T}.
{
$$\includegraphics[width=6cm]{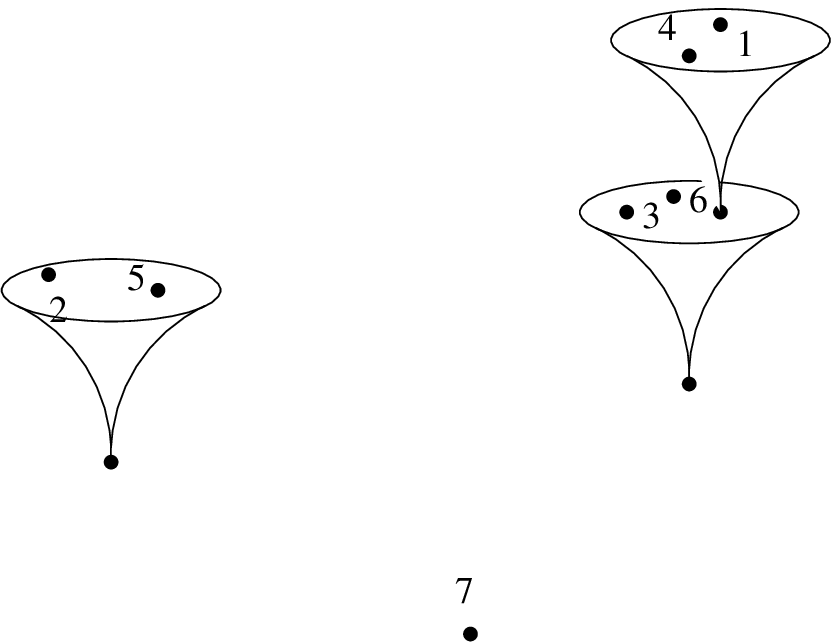}$$
}

\end{center}
\end{mydiagram}
\end{minipage}
\end{center}

%\begin{figure} \label{F:picDT}
%$$\includegraphics[width=6cm]{DT.eps} $$
%\caption{A tree $T$ and a representation of a point in $D_T(M)$.} 
%\end{figure} 
}

The following theorem is the main theorem of this section. 

\begin{theorem}\label{T:idstrata} 
$C_T(M)$ is diffeomorphic to $D_T(M)$. 
\end{theorem}

\begin{remark}
To intuitively understand $C_T(M)$ as part of the boundary of $C_n[M]$ one 
views an element of $IC_i(M)$ as a limit of a sequence in $C_i(M)$ which
approaches a point $(x, x, \ldots x)$ in the (thin)
diagonal of $M^i$.  Eventually, in such a sequence all the points in a
configuration would lie in a coordinate neighborhood of
$x$, which through the exponential map can be identified with $T_x M$,
and the limit is taken in that tangent space up to rescaling.  If
$i>2$, $IC_i(M)$ is itself not complete, so one allows these infinitesimal
configurations to degenerate as well, and this is how the situation is
pictured in Figure~\ref{F:pictDT}.  Because $T(TM) \cong \oplus_3 TM$, the
recursive structure of degenerating sub-configurations is not reflected
in the topology of $D_T(M)$.
\end{remark}

To establish this theorem we focus on the case in which $M$ is
Euclidean space $\R^m$, as $D_T(\R^m)$ admits a simple description.

\begin{definition}  \label{D:cntild}
Let $\wt{C}_n(\R^m)$ be the quotient of $C_n(\R^m)$ by $Sim_k$ acting 
diagonally, and let $q$ denote the quotient map.  
Choose coset
representatives to identify $\wt{C}_n(\R^m)$ with the subspace of 
$C_n(\R^m)$  of $(x_i)$ with $\Sigma_i x_i =  0$ and such 
that the maximum of the $d(x_i, \vec{0})$ is one. 
\end{definition} 

Because the tangent bundle of $\R^m$ is trivial, $IC_i(\R^m) \cong \R^m
\times \wt{C}_i(\R^m)$, and we have the following.

\begin{proposition}\label{P:DTRk}
$D_T(\R^m) = C_{\#v_0}(\R^m) \times  
\left(\wt{C}_{\#v}(\R^m)\right)^{V^i(T)}.$
\end{proposition}

Alternately, $D_T(\R^m)$ is the space in which each edge in $T$ is assigned
a point in $\R^m$, with coincident edges assigned distinct points, modulo
translation and scaling of coincident groups of edges.

Roughly speaking, the proof of \refT{idstrata} when $M = \R^m$ respects the
product decomposition of \refP{DTRk}.  We start by addressing the  stratum
associated to the tree $\includegraphics[width=0.4cm]{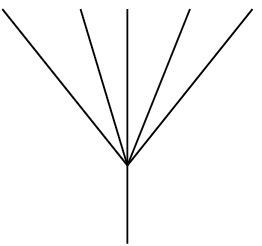}$ with a single 
internal vertex
connected to a univalent root.

\subsection{The auxilliary construction $C_n\{\R^m\}$}

\begin{definition}
Let $A_n\{M\} = (M)^{\n} \times (S^{m-1})^{C_2(\n)} \times 
(0, \infty)^{C_3(\n)}$, a subspace of $A_n[M]$.
Note that the image of $\alpha_n \colon C_n(M) \to A_n[M]$ lies
in  $A_n\{M\}$.  Let $C_n\{M\}$ be $cl_{A_n\{M\}}(im(\alpha_n))$.
\end{definition}

For our purposes, $C_n\{M\}$ will be useful as a subspace of $C_n[M]$ to
first understand, which we do for $M = \R^m$.

\begin{theorem}\label{T:firstcase}
$C_n\{ \R^m\}$ is diffeomorphic to 
$D_n\{\R^m\} = \R^m\times \wt{C}_n(\R^m) \times [0, \infty)$. 
\end{theorem}

As a manifold with boundary $C_n\{\R^m\}$ has two strata, namely 
$\R^m\times \wt{C}_n(\R^m) \times (0, \infty)$, which we will 
identify with $C_n(\R^m)$, and $\R^m\times \wt{C}_n(\R^m) \times 0$, 
the points added in this closure.
We will see that these correspond to 
$C_{\includegraphics[width=0.2cm]{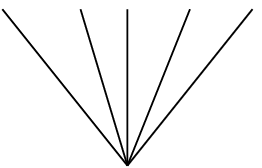}}(\R^m)$ and
$C_{\includegraphics[width=0.2cm]{rooted.eps}}(\R^m)$, respectively.
 
To prove Theorem~\ref{T:firstcase} we define a map 
$\nu \colon D_n\{\R^m\} \to A_n\{\R^m\}$ and show that it is a
homeomorphism onto $C_n\{ \R^m\}$. The map $\nu$ will essentially be an 
expansion from the point in $\R^m$ of the infinitesimal
configuration given by the point in $\wt{C}_n(\R^m) \subset
(\R^m)^{\und{n}}$.

\begin{definition}
\begin{enumerate}
\item Define $\eta \colon D_n\{\R^m\} \to (\R^m)^n$ by sending 
$x \times (y_i) \times t$ to $(x + t y_i)$.
\item Let $p$ denote the projection from $D_n\{\R^m\}$ onto
$\wt{C}_n(\R^m)$.
\item Let $\wt{\pi}_{ij}$ and $\wt{s}_{ijk}$ denote the maps on 
$\wt{C}_n(\R^m)$ which when composed with $q$ 
give the original $\pi_{ij}$ and $s_{ijk}$.
\item Finally, define $\nu \colon D_n\{\R^m\} \to A_n\{\R^m\}$ by 
$\eta \times (\wt{\pi}_{ij} \circ p) \times (\wt{s}_{ijk} \circ p)$.
\end{enumerate}
\end{definition}

When $t> 0$, the image of $\eta$ is in $C_n(\R^m)$, and moreover we have
the following.

\begin{proposition}\label{P:t>0} 
The map $\nu|_{t>0}$ coincides with 
$\alpha_n \circ \eta,$ a diffeomorphism from 
$\R^m\times \wt{C}_n(\R^m) \times (0, \infty)$ onto the image of 
$\alpha_n$.
\end{proposition}

\begin{proof}
For $t>0$, the map $\eta$ satisfies $\wt{\pi}_{ij} \circ p = \pi_{ij} \circ \eta$,
and similarly $\wt{s}_{ijk} \circ p  = s_{ijk} \circ \eta$, showing that
$\nu|_{t>0}$ coincides with $\alpha_n \circ \eta$.  

The inverse to $\nu|_{t>0}$ is the product of: the map which sends $(x_i)$ to its the
center of mass, the quotient map $q$ to $\wt{C}_n(\R^m)$, and the map
whose value is the greatest distance from one of the $x_i$ to the 
center of mass.   Both $\nu|_{t>0}$ and its inverse are clearly smooth.
\end{proof}

\begin{corollary}
$\nu|_{t=0}$ has image in  $C_n\{\R^m\}$.
\end{corollary}

We come to the heart of the matter, identifying $C_n\{\R^m\}$ 
as a closed subspace of $A_n\{\R^m\}$.  
We will apply this case repeatedly in analysis of $C_n[\R^m]$.
 
\begin{definition} 
Let $\wt{A}_n[\R^m] = (S^{m-1})^{C_2(\n)} \times I^{C_3(\n)}$, and let
$\wt{A}_n\{\R^m\} = (S^{m-1})^{C_2(\n)} \times (0,\infty)^{C_3(\n)}$. 
\end{definition} 

\begin{convention}
We extend multiplication to $[0,\infty]$by setting
$a \cdot \infty = \infty$ if $a \neq 0$ and $0 \cdot \infty = 1$.
\end{convention}

\begin{definition}
We say that vectors $\{ v_i \}$ are positively dependent if 
$\Sigma a_i v_i = 0$ for some $\{ a_i \}$ with all $a_i > 0$.
Similarly, $\{ v_i \}$ are non-negatively dependent if all $a_i \geq 0$.
\end{definition}

\begin{lemma}\label{L:lem1} 
The map 
$\iota_n = (\wt{\pi}_{ij}) \times (\wt{s}_{ijk}): \wt{C}_n(\R^m) \to 
\wt{A}_n[\R^m]$  is a diffeomorphism onto its image, which is closed as a
subspace of $\wt{A}_n\{\R^m\}$. 
\end{lemma} 

\begin{proof} 
Collinear configurations up to translation and scaling are cleary 
determined by their image under one  $\wt{\pi}_{ij}$ and the
$\wt{s}_{ijk}$.   
For non-collinear configurations, we may reconstruct ${\bf{x}} = (x_i)$
from the $u_{ij}= \wt{\pi}_{ij}({\bf{x}})$ and $d_{ijk} =
\wt{s}_{ijk}({\bf{x}})$ up to translation and scaling by for example
setting  $x_1 = \vec{0}$, $x_2 = u_{12}$ and then $x_i = d_{1i2} u_{1i}$
for any $i$.   These assigments of $x_i$ are smooth functions, so in fact
$\iota_n$ is a diffeomorphism onto its image.

For the sake of showing that the image of $\iota_n$ is closed, as well as
use in Section~\ref{S:simplicial}, we note that $d_{1i2}$ can be
determined  from the $u_{ij}$ by the law of sines.  If $\pm u_{ij}$, $\pm
u_{jk}$ and
$\pm u_{ik}$ are distinct, then 
$$\frac{|x_i - x_j|}{\sqrt{1 - (u_{ki} \cdot u_{kj})^2}} =
\frac{|x_j - x_k|}{\sqrt{1 - (u_{ij} \cdot u_{ik})^2}} =
\frac{|x_i - x_k|}{\sqrt{1 - (u_{ji} \cdot u_{jk})^2}}.
$$
Thus, in most cases $d_{1i2} = \sqrt{ \frac{{1 - (u_{2i} \cdot u_{21})^2}}
{{1 - (u_{i1} \cdot u_{i2})^2}}}$.  In general, as long as not all points
are collinear, the law of sines above can be used repeatedly to determine
all $d_{ijk}$ from the $u_{ij}$, which shows that when restricted to
non-collinear configurations, $(\wt{\pi}_{ij})$ itself is injective.

We identify the image of $\iota_n$ 
as the set of all points $(u_{ij}) \times (d_{ijk})$ 
which  satisfy the following conditions needed to consistently
define an inverse to $\iota_n$: 
\begin{enumerate} 
\item $u_{ij} = -u_{ji}$. \label{c0}
\item $u_{ij}$, $u_{jk}$ and $u_{ki}$ are positively dependent. 
\label{c1} 
\item If $\pm u_{ij}$,  $\pm u_{jk}$ and $\pm u_{ik}$ are distinct, then
$d_{ijk} = \sqrt{
\frac{{1 - (u_{ik} \cdot u_{jk})^2}} {{1 - (u_{ij} \cdot u_{jk})^2}}}$.
\label{c3} 
\item $d_{ijk}$ are non-zero and finite, and 
$$d_{ijk}  d_{ikj} = 1 = d_{ijk} d_{jki} d_{kij} =d_{ijk} d_{i \ell j} d_{i k \ell}.$$ 
\label{c2}  

\end{enumerate} 

We say a condition is closed if the subspace of points which satisfy it is
closed.  Note that Condition~\ref{c2} follows from Condition~\ref{c3}
when the latter applies.

Condition~\ref{c0} is clearly closed, and
Condition~\ref{c2} is a closed condition in
$\wt{A}_n\{\R^m\}$,  since we are already assuming that $d_{ijk} \in  (0,
\infty)$.  Condition~\ref{c3} says that on an open subspace of this image,
the $d_{ijk}$ are a function of the $u_{ij}$ and gives no restrictions
away from this subspace, and so is also a closed condition.  Considering 
Condition~\ref{c1}, it is a closed condition for $u_{ij}$,
$u_{jk}$ and  $u_{ik}$ to be dependent, but it is not usually closed to
be strictly positively dependent.   But the only dependence which can occur with
a coefficient of zero happens when $u_{ki} = -u_{jk}$ and $u_{ij} \neq \pm  u_{jk}$.
In this case $d_{ijk}$ would need to be $0$ by Condition~\ref{c3}, 
which cannot happen in $\wt{A}_n\{\R^m\}$.
So in fact Condition~\ref{c1} is closed within the points in
$\wt{A}_n\{\R^m\}$ satisfying Condition~\ref{c3}. 
\end{proof} 

Because $\nu|_{t=0}$ is the product of the diagonal map $\R^m\to
(\R^m)^{\n}$, which is a diffeomorphism onto its image,
 with $\iota_n$ we may deduce the following.

\begin{corollary}\label{C:0injective}
$\nu|_{t=0}$ is a diffeomorphism onto its image.
\end{corollary}

We may now finish analysis of $C_n\{\R^m\}$.

\begin{proof}[Proof of Theorem~\ref{T:firstcase}]
\refP{t>0} and \refC{0injective} combine to give that 
$\nu_T \colon D_T \{\R^m\} \to C_T \{ \R^m\}$ is injective.  We thus want 
to show that it is surjective and has a continuous inverse.

Consider the projection $p$ from $C_n\{\R^m\} \subset A_n\{\R^m\}$ to 
$(x_i) \in (\R^m)^{\n}$.  Over $C_n(\R^m)$ the image of
$\alpha_n$  is its graph, which is locally closed, so $p^{-1}(C_n(\R^m)) \cong
C_n(\R^m)$.  If $x_i = x_j$ but $x_i \neq x_k$ for
some $i,j,k$, continuity of $s_{ijk}$ would force $d_{ijk} = 0$, which is not
possible in $A_n\{\R^m\}$.  Thus no points in $C_n \{ \R^m\}$ lie over
such $(x_i)$. Over the diagonal of $(\R^m)^{\n}$ we know that $C_n \{\R^m
\}$  contains at least the image of $\nu|_{t=0}$.  But by \refL{lem1}, we
may  deduce that this image is closed in $A_n\{\R^m\}$ and thus accounts
for  all of $C_n\{\R^m\}$ over the diagonal.

We define an inverse to $\nu_T$ according to this decomposition over 
$(\R^m)^{\n}$. For a point in $C_n(\R^m)$, the inverse was given in
\refP{t>0}.  For points over the diagonal $(x_i = x)$ in
$(\R^m)^{\n}$,  the inverse is a product of the map which sends such a point
to $x \in \R^m$, 
$\iota_n^{-1}$, and the constant map whose image is $0 \in [0,\infty)$.
Smoothness of this inverse is straightforward and left to the reader.
\end{proof}

\subsection{Proof of \refT{idstrata} for $M = \R^m$}

Analysis of $C_T(\R^m)$ parallels that of $C_n\{\R^m\}$. A key 
construction is that of a map $\nu_T : N_T \to A_n[\R^m],$ 
where  $N_T\subset D_T(\R^m) \times [0,1)^{V^i(T)}$ is a 
chosen neighborhood of $D_T(\R^m) \times (t_v = 0)$. 
Although as mentioned before, $D_T(\R^m)$ is a subspace of $(\R^m)^{E(T)}$,
we emphasize the role of the vertices of $T$ in the definition of
$D_T(\R^m)$ by naming coordinates on
${{\bf x}} \in D_T(\R^m)$ as
${{\bf x}} = (x^v_e)$, where $v \in V(T)$ and $e \in e(V)$.
Recall that for each $v \neq v_0$ we consider $\wt{C}_{\#v}(\R^m)$ as
a subspace of $C_{\#v}(\R^m)$ in order to fix each $x^v_e$ as
an element of $\R^m$.

\begin{definition}\label{D:NT}
\begin{enumerate}
\item Let $N_T(\R^m)$ be the subset of $D_T(\R^m) \times [0,1)^{V^i(T)}$ 
of points ${{\bf x}} \times(t_v),$ where
${{\bf x}}$ can be any point in $D_T(\R^m)$, 
all $t_v < r({{\bf x}})$, defined by 
$$\frac{r({\bf x})}{(1-r({\bf x}))} = 
\frac{1}{3}{\text{min}} \{d(x^v_e, x^v_{e'} ) \}, \;\; {\text{where}} \;\;
{v \in V(T), \;\; e, e' \in E(v)}.$$

\item By convention, set $t_{v_0} = 1$.  Let $s_w : N_T \to [0,1)$ send  
${{\bf x}} \times(t_v)$ to the product of $t_v$ for $v$ in the root path of $w$.

\item
For any vertex $v$ of an $f$-tree $T$ define $y_v \colon N_T(\R^m) \to
\R^m$ inductively by setting $y_{v_0} = 0$ and 
$y_v({\bf{x}}) = s_{w} x^{w}_e + y_{w}({\bf{x}})$,  where $e$ is the
edge for which $v$ is terminal and $w$ is the initial vertex of $e$.
Define
$\eta_T : N_T \to (\R^m)^{l(T)}$ to be $(y_\ell)^{l(T)}$.

\end{enumerate}
%Given an $f$-tree $T$ and a leaf $\lambda$ of $T$ let $r(\lambda)$ be the 
%root edge over which that $\lambda$ lies, let 
%$v(\lambda)$ be the terminal vertex of $r(\lambda)$, and let $o(\lambda)$ be
%the place of leaf $\lambda$ in the ordering of leaves of $T_{v(\lambda)}$ (if
%$T_{v(\lambda)}$ is  non-trivial, that is, $\lambda \neq v(\lambda)$).
%
%Define $\eta_T : N_T \to (\R^m)^{l(T)}$ inductively as
%$\eta_T({\bf x}, (t_v)) = (y_i)$ with
%$$y_i = 
%\begin{cases}
%x^{v_0}_{r(\lambda)}  & {\text{if}} \; \lambda = v(\lambda)
%\\ x^{v_0}_{r(\lambda)} + t_{v(\lambda)} \cdot x_{o(\lambda)}
%      \left(\eta_{T_{v(\lambda)}} \circ \pi_{v(\lambda)} ({\bf x}, (t_v)) \right)
%            &  {\text{otherwise}}.
%\end{cases}
%$$
%Here we recall that $x_{o(\lambda)}$ is the $o(\lambda)$th coordinate function of 
%$(\R^m)^{l(T_{v(\lambda)})}$.
\end{definition}

\skiphome{

\begin{center}
\begin{minipage}{10cm}
\begin{mydiagram}\label{F:etaT}\begin{center}

{$\eta_T$  of the point from Figure~\ref{F:pictDT} (and some $t_v > 0$)
with all $y_v$ indicated.} 
\psfrag{A}{$y_{v_0}$}
\psfrag{B}{$y_{v_1}$}
\psfrag{C}{$y_{v_2}$}
\psfrag{D}{$y_{v_3}$}
\psfrag{E}{$y_1$}
\psfrag{F}{$y_4$}
\psfrag{G}{$y_7$}
\psfrag{H}{$y_6$}
\psfrag{I}{$y_3$}
\psfrag{J}{$y_2$}
\psfrag{K}{$y_5$}
$$\includegraphics[width=6cm]{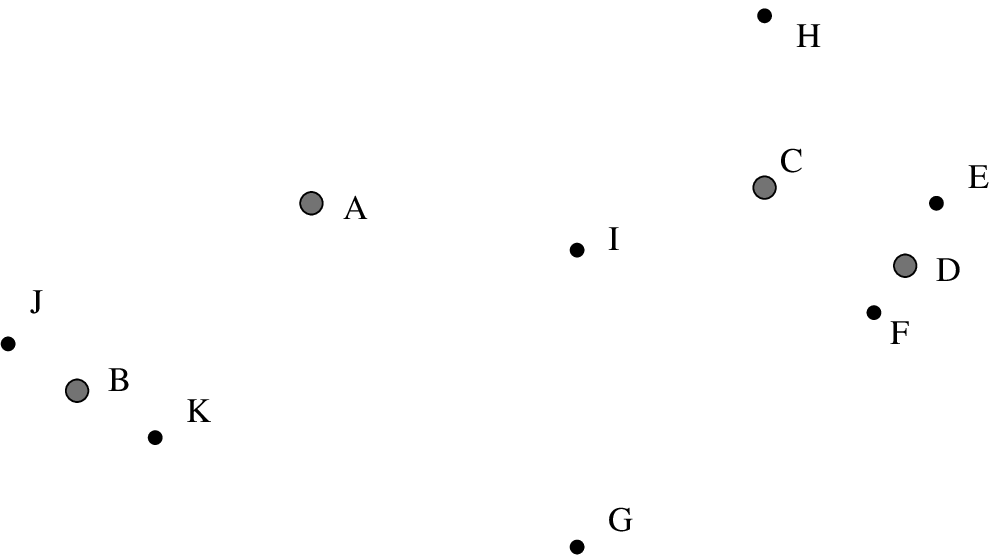} $$

\end{center}
\end{mydiagram}
\end{minipage}
\end{center}

}

See Figure~\ref{F:etaT} for an illustration of this construction.
The most basic case is when $T = \includegraphics[width=0.4cm]{terminal.eps}$ 
the terminal object of $\Psi_{\und{n}}$, in which case
$N_{\includegraphics[width=0.2cm]{terminal.eps}}(\R^m) = 
D_{\includegraphics[width=0.2cm]{terminal.eps}}(\R^m) = C_n(\R^m)$  and
$\eta_{\includegraphics[width=0.2cm]{terminal.eps}}$ is the canonical inclusion in
$(\R^m)^{\und{n}}$.

\begin{definition}
\begin{enumerate}
\item Given a vertex $w$ of $T$, let $T_w$ be
the $f$-tree consisting of all vertices and edges over $w$, where $w$
serves as the root of $T_w$ and the leaves over $w$ are re-labelled 
consistent with the order of their labels in $T$.   

\item Let $\rho_w \colon N_T(\R^m) \to N_{T_w}(\R^m)$ be the
projection onto factors indexed by vertices in $T_w$, but with $t_w$ set to 
one and the projection onto $\wt{C}_{\#w}(\R^m)$ composed with  the canonical
inclusion to $C_{\#w}(\R^m)$.

\item Let $f_{ij}$ be the composite $\pi_{i'j'} \circ
\eta_{T_v} \circ \rho_v$, where $v$ is the join of the leaves labelled
$i$ and $j$, and $i'$ and $j'$ are the labels of the corresponding
leaves of $T_v$.  Similarly, let $g_{ijk}$ be the composite $s_{i'j'k'}
\circ \eta_{T_w} \circ \rho_w$ where $w$ is the join of leaves $i$, $j$ and $k$.

\item 
Define $\nu_T : N_T \to A_n[\R^m]$ to be the product $\eta_T \times (f_{ij})
\times (g_{ijk})$.   
Let $\nu^0_T$ be the restriction of $\nu_T$ to $D_T \times (0)^n \subset N$.
\end{enumerate}
\end{definition}

\begin{proposition}\label{P2}
The image of $D_T(\R^m)$ under $\nu^0_T$ lies in $C_T(\R^m)$.
\end{proposition}

\begin{proof}
First note that $\eta_T|_{(t_i > 0)}$ has image in $C_n(\R^m)$. 
Moreover,  if all $t_i > 0$, $f_{ij}$
coincides with $\pi_{ij} \circ  \eta_T$ and similarly 
$g_{ijk} = s_{ijk} \circ \eta_T$.  Thus, the image of 
$\nu_T|_{(t_i > 0)}$ lies in $\alpha_n(C_n(M))$, which implies that 
all of the image of $\nu_T$, and in particular that of $\nu^0_T$, 
lies in $C_n[M]$.

If $v$ is the join of leaves $i,j$ and $k$ and we set
$(y_i) = \nu^0_{T_v} \circ \rho_v ({{\bf x}}, (t_v))$, then
$y_{i'} = y_{j'}$  and thus  $d_{ijk} = s_{ijk}((y_i)) = 0$ if and only if the join of 
leaves $i$ and $j$ is some vertex which lies (strictly) over $v$.   
Thus, the exclusion relation for $\nu_T^0({\bf x}, (t_v))$ 
as an element of $C_n[M]$ is the exclusion relation associated to $T$.
\end{proof}

The simplest way to see that $\nu^0_T$ is a homeomorphism onto
$C_T(\R^m)$ is to decompose it as a product and use our analysis of 
$C_n\{\R^m\}$ to help define an inverse.

\begin{definition}
Let $A_T[M] \subset A_n[\R^m]$ be the subset of points
$(x_i) \times (u_{ij})\times (d_{ijk})$ such that 
\begin{enumerate}
\item If the join of leaves $i$ and $j$ is not the root vertex, then
$x_i = x_j$.
\item If $(i,j),k$ is in the exclusion relation ${\mathcal{E}x}(T)$,
then $d_{ijk} = 0$, $d_{ikj} = \infty$, $d_{kij} = 1$, 
and $u_{ik} = u_{jk}$.
\end{enumerate}
Let $A_T\{M\}$ be the subspace of $A_T[M]$ for which 
if there are no exclusions among $i$, $j$ and $k$, then  $d_{ijk}$
is non-zero and finite.
\end{definition}

We claim that $C_T(\R^m) = C_n[\R^m]  \cap A_T\{\R^m\}$.
The relations between
the $x_i$, $u_{ij}$ and $d_{ijk}$ which hold on $C_n(\R^m)$ also hold
on $C_T(\R^m)$ by continuity. 
Therefore, the defining conditions of $A_T\{\R^m\}$ when restricted
to its intersection with $C_n[\R^m]$
will follow from the conditions 
$d_{ijk} =0$ when $(i,j)k \in {\mathcal{E}x}(T)$, which in turn are the
only defining conditions for $C_T[\R^m]$.

Thus the image of $\nu^0_T$ lies in $A_T\{\R^m\}$.
By accounting for diagonal subspaces and reordering terms, we will decompose
$\nu^0_T$ as a product of maps in order to define its inverse.
We first set some notation.

\begin{definition}\label{D:funct}
Given a map of sets $\sigma : R \to S$ let
$p^X_{\sigma}$, or just $p_{\sigma}$,
denote the map from $X^S$ to $X^R$ which sends 
$(x_i)_{i \in S}$ to $(x_{\sigma(j)})_{j \in R}$.
\end{definition}

\begin{definition}

\begin{enumerate}
\item Given a tree $T$ choose $\sigma_0 : \und{\#v_0} \to \und{n}$ to be an inclusion
of sets such that each point in the image labels a leaf which lies over a
distinct root edge of $T$.  

\item Similarly, choose $\sigma_v : \und{\#v} \to
\und{n}$ to be an inclusion whose image labels leaves which lie over distinct
edges for which $v$ is initial.  

\item Let $p_{v_0} : A_n[\R^m] \to A_{\#v_0}[\R^m]$ be the projection
$p_{\sigma_0} \times p_{C_2(\sigma_0)} \times p_{C_3(\sigma_0)}$.
\item Similarly, let $p_v : A_n[\R^m] \to \wt{A}_{\#v}[\R^m]$ be the 
projection $* \times p_{C_2(\sigma_v)} \times p_{C_3(\sigma_v)}$.

\item Let $p_T$ be the product $(p_v)^{V(T)}$.
\end{enumerate}
\end{definition}

For example, with $T$ as in Figure~\ref{F:T},
the image of $\sigma_0$ could be $\{ 5, 7, 4\}$ and  of $\sigma_{v_1}$ 
could be $\{2, 5\}$.

\begin{proposition}\label{P1}
For any choice of $\sigma_v$, the projection
$p_T$ restricted to
$A_T[\R^m]$ is a diffeomorphism onto $A_{\#v_0}[\R^m] \times \left(
\wt{A}_{\#v}(\R^m)\right)^{V^i(T)}$, splitting the inclusion of
$A_T[\R^m]$ in $A_n[\R^m]$.   Moreover, composed with this diffeomorphism, 
$\nu^0_T$ is the product $\alpha_{\#v_0} \times (\iota_{\#v})$.
Analogous results hold for $A_T\{\R^m\}$.
\end{proposition}

We leave the proof of this proposition, which is essentially
unraveling definitions, to the reader. 
We will now define the inverse to $\nu_T^0$ one vertex at a time.
For $v \in V^i(T)$, consider $p_v(y) \in \wt{A}_{\#v}\{\R^m\}$, which by \refL{projclose}
lies in the closure of the image of $p_v|_{C_n(\R^m)}$.  The image of  
$p_v|_{C_n(\R^m))}$ coincides with the image of
$\iota_{\#v}$, and by \refL{lem1}  the image of $\iota_{\#v}$ is
already closed in $\wt{A}_{\#v}\{\R^m\}$.   Moreover, $\iota_{\#v}$  is a
diffeomorphism onto its image, so we may define the following.

\begin{definition}\label{D:phiT}
\begin{enumerate}
\item For $v \in V^i(T)$, let $\phi_v \colon C_T(\R^m) \to \wt{C}_{\#v}(\R^m) = 
\iota_{\#v}^{-1} \circ p_v$.
\item For $v_0$, note that if $y \in C_T(\R^m)$, then $p_{v_0}(y)$ lies in the image of 
$\alpha_{\#v_0}$.  Define $\phi_{v_0} = \alpha_{\#v_0}^{-1} \circ p_{v_0}$.
\item Let $\phi_T  = (\phi_v)^{v \in V(T)} :  C_T(\R^m) \to C_{\#v_0}(\R^m) \times
\left( C_{\#v}(\R^m) \right)^{V^i(T)}$.
\end{enumerate}
\end{definition}

In other words, $\phi_T$ is the composite:
\begin{multline}
C_T(\R^m) \subset A_T \{ \R^m\} \overset{p_v}{\to} 
    A_{\#v_0}\{\R^m\} \times \left(\wt{A}_{\#v}\{\R^m\}  \right)^{V^i(T)} \\
\overset{\alpha_{\#v_0}^{-1} \times (\iota_{\#v}^{-1})^{V^i(T)})}{\longrightarrow} 
C_{\#v_0}(\R^m) \times \left(\wt{C}_{\#v}(\R^m) \right)^{V^i(T)} = D_T(\R^m).
\end{multline}

\begin{proof}[Proof of \refT{idstrata} for $M = \R^m$]
By Proposition~\ref{P2}, $\nu^0_T$ sends $D_T(\R^m)$ to 
$C_T(\R^m) \subset A_n[\R^m]$.
\refD{phiT} constructs $\phi_T : C_T(\R^m)  \to D_T(\R^m)$. 
By construction, and appeal to \refL{lem1}, they
are inverse to one another.  We also need to check that  $\nu_T^0$ and
$(\phi_v)$ are smooth, which follows by checking that their component functions
only involve addition, projection and  $\iota_n^{-1}$ which we know
is smooth from \refL{lem1}.
\end{proof}

\subsection{Proof of \refT{idstrata} for general $M$}

To establish \refT{idstrata} for general $M$ we first identify $D_T(M)$ as a 
subspace of $D_T(\R^m)$, and then we will  make use of the established diffeomorphism 
between $D_T(\R^m)$ and $C_T(\R^m)$.   To set notation, let $\epsilon$ be the given 
embeddding
of $M$ in $\R^m$.
 
%For $M$, a submanifold of $\R^m$,   
%$C_T(M) = C_n[M] \cap C_T(\R^m)$ in $A_n[M] \subset A_n[\R^m]$ almost by
%definition.  Set coordinates for points in  $D_T(\R^m)$ as $(x^{v_0}_e)_{e \in
%E_0} \times ({\bf{x}}^v)_{v \in V^i(T)}$ where $(x^{v_0}_e) \in C_{\#v_0}(\R^m)$ and
%${\bf{x}}^v \in \wt{C}_{\#v}(\R^m)$.

\begin{proposition}\label{P:idDTM}
The subspace $ID_n(M)$
of $M \times \wt{C}_n(\R^m)$ consisting of all $(m,{\bf x})$ such that all
$\pi_{ij}({\bf x})$ are in $T_m M \subset \R^m$ is diffeomorphic, 
as a bundle over $M$, to $IC_n(M)$.
Through these diffeomorphisms, $D_T(M)$ is a subspace of $D_T(\R^m)$.
\end{proposition}

\begin{proof}[Proof sketch.]
The first statement follows from the
standard identification of $TM$ as a sub-bundle of $T\R^m|_M$.
The second statement follows from the first statement 
and Definition~\ref{D:DT} of $D_T(M)$.
\end{proof}

From now on, we identify $D_T(M)$ with this subspace of $D_T(\R^m)$.

\begin{proposition}\label{P:P3}
$C_T(M) \subseteq \nu^0_T\left(D_T(M)\right)$.
\end{proposition}

\begin{proof}
Since $C_T(M)$ is already a subspace of $C_T[\R^m]$, we just need to check
that its points satisfy the condition of \refP{idDTM}.  Looking at
$(x_i) \times (u_{ij}) \times (d_{ijk}) \in C_n(M)$ inside 
$A_n[M]$  we see that the $u_{ij}$ are vectors which are secant to $M$.
Thus, in the closure, if $x_i = x_j$, then $u_{ij}$ is tangent to
$M$ at $x_i$. 
\end{proof}

To prove the converse to this proposition, we show that $\nu^0_T(D_T(M))$
lies in $C_T(M)$ by modifying the maps $\eta_T$ and $\nu_T$ so that the
image of the latter is in the image of $\alpha_n$.   The easily remedied
defect of $\eta_T$ is that it maps to the image of the tangent bundle of $M$ in $\R^m$, 
and not to $M$ itself.  

\begin{definition}
\begin{enumerate}
\item Let $N_T(M) \subset N_T(\R^m)$ be the subspace of $({\bf x}, (t_v))$ with ${\bf x} \in D_T(M)$.  
\item Let $\eta^*_T: N_T(M) \to (\R^m)^n \times (\R^m)^n = T(\R^m)^n$ send 
$({\bf x}, (t_v))$ to $\eta_T({\bf x}, (0)) \times \eta_T({\bf x}, (t_v))$.
\end{enumerate}
\end{definition}

\begin{proposition}
The image of $\eta^*_T$ lies in $T\epsilon^n : TM^n \subset T(\R^m)^n$.
\end{proposition}

\begin{proof}
$\eta_T({\bf x}, (t_v))$ is defined by adding vectors which by \refP{idDTM} are tangent
to $M$ to the coordinates of $\eta_T({\bf x}, (t_v))$, which are in $M$.
\end{proof}

We map to $M^{\und{n}}$ by composing with the exponential map $Exp(M)$.
For each ${\bf x} \in D_T(M)$ let
$U_{\bf{x}}$ be a neighborhood of ${\bf x} \times 0$ in $N_T(M)$ such that
the exponential map $Exp({M^{\und{n}}})$ is injective on 
$\eta^*_T(U_{\bf{x}})$.

\begin{definition}
\begin{enumerate}
\item Let $\eta^{M, {\bf{x}}}_T \colon U_{\bf{x}} \to M^{\und{n}}$ 
be the composite  $Exp({M^{\und{n}}}) \circ (T\epsilon^n)^{-1} \circ \eta_T^*$.
\item Define $f_{ij}^{M, {\bf x}}$ 
by letting $(z_i)$ denote $\eta_T^{M, {\bf{x}}}({\bf{y}}, (t_v))$
and setting $f_{ij}^{M, {\bf x}} = \pi_{ij} \circ \eta_T^{M, {\bf x}}$ if $z_ i \neq z_j$ or
$DExp \circ f_{ij}$, where $DExp$ is the derivative of the exponential map 
at $z_i \in TM$ and the composite is well defined since $TM \subset T\R^m$.

\item Define $\nu^{M,{\bf{x}}} _T :  U_{{\bf x}} \to A_n[M]$ as the product  
$\eta_T^{M, {\bf{x}}} \times (f^M_{ij}) \times 
                          (s_{ijk}  \circ \eta_T^{M, {\bf{x}}})$.  
\end{enumerate}
\end{definition}

By construction, the image of $\nu^{M,{\bf{x}}}_T|_{(t_i > 0)}$ in $A_n[M]$
lies in the image of $\alpha_n$. 
On $D_T(M) \cap U_{\bf{x}}$ 
the map  $\nu^{M,{\bf{x}}}_T$ coincides with $\nu^0_T$
establishing  that $\nu_T^0(D_T(M)) \subset C_T(M)$.  Along with
\refP{P3} and the fact that $\nu_T^0$ and its inverse are smooth, this
completes the proof of \refT{idstrata}.

\section{First properties}

Having proved \refT{idstrata} we derive some first consequences from both
the theorem and the arguments of its proof.  

\subsection{Characterization in $A_n[M]$ and standard projections}

To map from $C_n[M]$, as we have defined it, one may simply restrict maps from
$A_n[M]$.  To map into $C_n[M]$ is more difficult, but the following theorem 
gives conditions to verify that some point in $A_n[M]$ lies in $C_n[M]$.

\begin{theorem}\label{T:cn[sub}
$C_n[\R^m]$ is the subspace of $A_n[\R^m]$ of points $(x_i) \times (u_{ij}) 
\times (d_{ijk})$ such that
\begin{enumerate}
\item If $x_{i} \neq x_k$, then $u_{ij} = \frac{x_i - x_k}{||x_i - x_k||}$ and 
$d_{ijk} = \frac{d(x_i, x_j)}{d(x_i, x_k)}$. \label{d1}
\item  If $\pm u_{ij}, \pm u_{jk},$ and $\pm u_{ik}$ are all distinct, then 
$d_{ijk}= \sqrt{\frac{{1 - (u_{ki} \cdot u_{kj})^2}}{{1 - (u_{ji} \cdot u_{jk})^2}}}.$
Otherwise, if $u_{ik} = u_{jk} \neq u_{ij}$, then $d_{ijk} = 0$. \label{d2}
\item $u_{ij} = -u_{ji}$, and $u_{ij}$, $u_{jk}$, and $u_{ki}$ are non-negatively
dependent. \label{d3}
\item $d_{ijk} d_{ikj}$, $d_{ijk} d_{ikl} d_{ilj}$ and $d_{ijk} d_{jki} d_{kij}$ 
are all equal to one. \label{d4}
\end{enumerate}
Moreover, $C_n[M]$ is the subspace of $C_n[\R^m]$ where all $x_i \in M$ and 
if $x_i = x_j$, then $u_{ij}$ is tangent to $M$ at $x_i$.
\end{theorem}

\begin{proof}
It is simple to check that $C_n[\R^m]$ satisfies all of the properties listed. 
In most cases, the properties are given by equalities which hold on $C_n(\R^m)$
and thus $C_n[\R^m]$ by continuity.  We noted in \refP{P3} that if $x_i =
x_j$ in $C_n[M]$, then $u_{ij}$ is tangent to $M$.  

Conversely, we can start with a point
$x$ which satisfies these properties, and 
Condition~\ref{d4} allows us to define $T(x)$ as in Definition~\ref{D:Tx}. 
We can then either mimmic the construction of $\eta_{T(x)}$ to find points in
the image of $\alpha_n$ nearby showing that $x \in C_n[M]$,
or go through the arguments of
Section~\ref{S:strata}, in particular the proof of \refL{lem1}, 
to find an element of $D_T(M)$ which maps to
$x$ under $\nu^T_0$.  The latter argument proceeds by showing that $x$ lies in
$A_{T(x)}[\R^m]$, as we may use the contrapositives to Conditions~\ref{d1} and
\ref{d2} along with \ref{d3} to show that if $d_{ijk} = 0$, then $x_i = x_j$
and $u_{ik} = u_{jk}$.  Then, $p_{v_0}(x)$ lies in the image of
$\alpha_{\#v_0}$ essentially by Condition~\ref{d1}.  Next, $p_v(x)$ lies in the
image of $\iota_{\#v}$ by Conditions~\ref{d2},~\ref{d3} and~\ref{d4}, as these
conditions coincide with those given for the image of $\iota$ in \refL{lem1}. 
We apply the product $(\phi_v)$ to get a point in $D_T(\R^m)$ which
maps to $x$ under $\nu_T^0$.
\end{proof}

We next turn our attention to the standard projection maps.  

\begin{theorem}
By restricting the projection of $A_n[M]$ onto $(M)^n$ to $C_n[M]$,
we obtain a projection map $p$ which is onto, which extends the inclusion
$\iota$ of $C_n(M)$ in $(M)^n$ and for which every point in $C_n(M)$ has 
only one pre-image.
\end{theorem}

\begin{proof}
The fact that $p$ is onto can be seen through composing $p$ with the maps
$\nu_T^0$.  It is immediate from definitions that $p$ extends $\iota$.  
Finally, by our characterization in
\refT{cn[sub}, in particular Condition~\ref{d1}, any point in $C_n[M]$ which
projects to $C_n(M)$ will be in the image of $\alpha_n$.
\end{proof}

When $M = \R^m$, it is meaningful to project onto other factors of
$A_n[\R^m]$ to get similar extensions.

\begin{theorem}
The maps $\pi_{ij}$ and $s_{ijk}$ extend to maps from $C_n[\R^m]$.
Moreover, the extension $\pi_{ij}$ is an open map.
\end{theorem}

\begin{proof}
The only statement which is not immediate is that $\pi_{ij}$ is 
an open map.  We check this on each stratum, using the
identification of $C_T(\R^m)$ as a product.  When
$\pi_{ij}$ is restricted to $C_T(\R^m)$ it factors as $p_v$,
where $v$ is the join of leaves $i$ and $j$, composed with 
some $\wt{\pi}_{i'j'}$, each of which is an open map.
\end{proof}

\subsection{Manifold structure, codimensions of strata, functoriality for embeddings, and equivariance}

\begin{theorem}\label{T:mfld}
$C_n[M]$ is a manifold with corners for which the 
$\nu^{M,{\bf{x}}}_T$ may serve as charts.
\end{theorem}

\begin{proof}
The domains of $\nu^{M,{\bf{x}}}_T$ are manifolds with corners,
so it suffices to check that these maps are diffeomorphisms onto their images in
$A_n[M]$, which is itself a manifold with corners.  We have already noted that
$\nu^{M,{\bf{x}}}_T$ are smooth on their domain, as they are defined using addition
in $\R^m$, projection maps, and the exponential map.  Moreover, they may be extended
using the same formulas to values of
$t_v < 0$, as needed for smoothness with corners.  

For $M = \R^m$, the inverse to $\nu_T$ is relatively straightforward to define.  
Given ${\bf {x}} = (x_i) \times (u_{ij}) \times (d_{ijk}) \in A_n[\R^m]$, first
recursively set $y_v$ to be the average of $y_w$, where $w$ are terminal
vertices for edges coincident at $v$, starting with $y_l = x_i$ when $l$
is the leaf labelled by $i$.  We let $(y_v)$, as $v$ ranges over terminal
vertices for root edges of $T$, define ${\bf{x}}_{v_0} \in C_{\#v_0}(\R^m)$.

Along the same lines, for each vertex $v$ first define a point in
$(\R^m)^{l(T_v)}/Sim_k$, as in the definition of
$\iota_n^{-1}$, by setting some $x_i = 0$, some $x_k = u_{ik}$ and the rest of the
$x_j$ as $d_{ijk} x_{ij}$.  Recursively set $x_w$ to be the average of
$x_u$ for $u$ directly over $w$ (which is well defined up to translation and
scaling) and let $(x_w)$ as $w$ ranges over vertices directly over
$v$ define ${\bf{x}}_v \in \wt{C}_{\#v}(\R^m)$.  

Finally, to compute $t_w$ we look within the construction above of
${\bf{x}}_v$ for the vertex
$v$ over which $w$ sits directly.  Let $d_w$ be the greatest distance from one of
the $x_u$, for $u$ over $w$, to $x_w$, and define $d_v$ similarly.  Set
$t_w$ to be $d_w/d_v$.  

The map which sends ${\bf{x}}$ as above to
$({\bf{x}}_v) \times (t_v)$ is the inverse to $\nu_T$, and it is smooth, defined
by averaging and greatest distance functions. The construction for general $M$
works similarly, by first composing with the inverse to the exponential map.  We
leave its construction to the reader.
\end{proof}

Since a manifold with corners is a topological manifold with boundary, and
a topological manifold with boundary is homotopy equivalent to its interiors, 
we get the following.

\begin{corollary}
The inclusion of $C_n(M)$ into $C_n[M]$ is a homotopy equivalnce.
\end{corollary}

An essential piece of data for a manifold with corners are the dimensions of the 
strata.
Dimension counting for $D_T(M)$ leads to the following.

\begin{proposition}
The codimension of $C_T(M)$ is $\# V^i(T)$.
\end{proposition}

Contrast this with the image of the projection of $C_T(M)$ in $M^{\und{n}}$,
which has codimension equal to $k \cdot dim(M)$, where $k$ is the sum over all
root edges $e$ of $n_e - 1$ where $n_e$ is the number of leaves over $e$.

Next, we have the following long-promised result.

\begin{theorem}\label{T:diff}
Up to diffeomorphism, $C_n[M]$ is independent of the embedding of $M$
in $\R^m$.
\end{theorem}

\begin{proof}
First note that the definitions of
$D_T(M)$ and $N_T(M)$ do not use the embedding of $M$ in $\R^m$. 
Let $f$ and $g$ be two embeddings of $M$ in $\R^m$, and
let $\nu^{f,{\bf{\und{x}}}}_T$ and $\nu^{g,{\bf{\und{x}}}}_T$ denote the
respective versions of $\nu^{M,{\bf{\und{x}}}}_T$.  Then the 
$\nu^{f,{\bf{\und{x}}}}_T \circ (\nu^{g,{\bf{\und{x}}}}_T)^{-1}$
compatibly define a diffeomorphism between the two versions of $C_n[M]$.
\end{proof}

Since the exponential maps from $T(M)^n$ to
$(M)^n$ are independent of the embedding of $M$ in $\R^k$, so are
$\nu^{M,{\bf{\und{x}}}}_T$. 
Thus, we could use the $\nu^{M,{\bf{\und{x}}}}_T$ to
topologize the union of the $C_T(M)$ without reference to $A_n[M]$.
Yet another approach would be to first develop
$C_n[\R^m]$ and then use a diffeomorphism result for these 
to patch $C_n[M]$ together from $C_m[U_i]$ for $m \leq n$, where $U_i$ is a
system of charts for $M$.  

\begin{corollary}
$C_n[-]$ is functorial in that an embedding $f \colon M \to N$ 
induces an embedding of manifolds with corners $C_n[f] \colon C_n[M] \to
C_n[N]$ which respects the stratifications.  Moreover, $C_T(M)$
is mapped to $C_T(N)$ by $Df$ on each factor of $IC_i(M)$.
\end{corollary}

\begin{proof}
Since we are free to choose the embedding of $M$ in $\R^m$ to define
$C_n[M]$, we may simply compose the chosen embedding of $N$ in $\R^m$ 
with $f$, giving immediately that $C_n[M]$ is a subspace of $C_n[N]$.
Moreover, by definition of the stratification according to conditions
of $d_{ijk} = 0$, $C_T(M)$ is a subspace of $C_T(M)$.
The fact that $C_n[M]$ is embedded as a submanifold with corners
is readily checked on each stratum, using the fact that $IC_i(M)$
is a submanifold of $IC_i(N)$ through $Df$.
\end{proof}

An alternate notation for $C_n[f]$ is $ev_n(f)$ as it extends the evaluation
maps on $C_n(M)$ and $M^{\und{n}}$.

\begin{corollary}\label{C:diff}
The group of diffeomorphisms of $M$ acts on $C_n[M]$, extending and lifting its
actions on $C_n(M)$ and $M^{\und{n}}$.
\end{corollary}

The construction of $C_n[M]$ is also
compatible with the free symmetric group action $C_n(M)$.

\begin{theorem}
The $\Sigma_n$ action on $C_n(M)$ extends to one on $C_n[M]$, which
is free and permutes the strata by diffeomorphisms according to the
$\Sigma_n$ action on $\Psi_{\und{n}}$.  Thus, the quotient $C_n[M]/ \Sigma_n$ is
itself a  manifold with corners whose category of strata is isomorphic to
$\wt{\Psi}_n$, the category of unlabelled $f$-trees.
\end{theorem}

\begin{proof}
The $\Sigma_n$ action on $C_n[M]$ may in fact be defined as the
restriction of the action on $A_n[M]$ given by permutation
of indices.

The fact that this action is free follows either from a
stratum-by-stratum analysis or, more directly, from the fact that if 
$\sigma$ is a permutation with a cycle $(i_1, \ldots, i_k)$  with $k>1$ 
and  if $u_{i_1i_2} = u_{i_2 i_3} = \cdots = u_{i_{k-1}i_{k}}= u$, then
$u_{i_1 i_k} = u$ as well by Condition~\ref{d3}  of \refT{cn[sub}. This implies 
that $u_{i_k i_1} = -u \neq u$, which means that $\sigma$ cannot fix
a point in $C_n[M]$ unless it is the identity.  

Finally, note that the coordinate charts $\nu_T^{M, \bf{\und{x}}}$ 
commute with permutation of indices, so that $\sigma C_T(M) = 
C_{\sigma T}(M)$ through a diffeomorphism, giving rise to a manifold
structure on the quotient.
\end{proof}

\subsection{The closures of strata}

We will now see that the passage from the stratum $C_T(M)$ to its closure, which by
\refT{strata1} consists of the union of $C_S(M)$ for $S$ with a morphism 
to $T$, is similar to the construction of  $C_n[M]$ itself.  Recall \refD{cntild} of
$\wt{C}_n(\R^m)$.

\begin{definition}
Let $\wt{C}_n[\R^m]$ be defined as the closure of $\wt{C}_n(\R^m)$ 
in $\wt{A}_n[\R^m]$.
\end{definition}

Because $\wt{A}_n[\R^m]$ is compact, so is $\wt{C}_n[\R^m]$.
We give an alternate construction of this space as follows.
Extend the action of $Sim_k$ on $(\R^m)^{\und{n}}$ to 
$A_n[\R^m]$ by acting trivially on the factors of $S^{m-1}$ and $I$.
This action preserves the image of $\alpha_n$ and so passes to an action
on $C_n[\R^m]$.  This is a special case of \refC{diff}.  Let 
$A_n[\R^m]/\sim$ and $C_n[\R^m]/\sim$ denote the quotients by these actions.

\begin{lemma}\label{L:sim}
$\wt{C}_n[\R^m]$ is diffeomorphic to $C_n[\R^m]/\sim$.
\end{lemma}

\begin{proof}
First note that $A_n[\R^m]/\sim$ is compact, and thus so is $C_n[\R^m]/\sim$. 
The projection map from $A_n[\R^m]/\sim$ to $\wt{A}_n[\R^m]$ thus 
sends $C_n[\R^m]/\sim$ onto $\wt{C}_n[\R^m]$ by \refL{projclose}.
In the other direction, we may essentially use the maps $\iota_k^{-1}$ to define
an inverse to this projection, by reconstructing a point in 
$(\R^m)^{\und{n}}$ up to translation and scaling 
from its images under $\pi_{ij}$ and $s_{ijk}$.
\end{proof}

We may define a stratification of $\wt{C}_n[\R^m]$ labelled
by trees in the same fashion as for $C_n[\R^m]$,
and the strata have a more uniform description than that of $C_T(\R^m)$.

\begin{corollary}\label{C:wtCT}
$\wt{C}_T(\R^m)$ is diffeomorphic to 
$\left(\wt{C}_{\#v}(\R^m)\right)^{V(T)}$.
\end{corollary}

\begin{proof}
We cite \refL{sim} and check that $Sim_k$ is acting on each $C_T(\R^m)$ 
non-trivially
only on the factor of
$C_{\#v_0}(\R^m)$, and doing so there by its standard diagonal action.
\end{proof}

Other results for $C_n[\R^m]$ have similar analogues for
$\wt{C}_n[\R^m]$, which we will not state in general.  One of note is that
its category of strata is isomorphic to $\wt{\Psi}_n$, the category of
trees with a trunk.

\begin{definition}
\begin{enumerate}
\item Define $IC_n[M]$ as a fiber bundle over $M$ with fiber
$\wt{C}_n[\R^m]$  built from $TM$  by taking the same system of charts
but choosing coordinate transformations $\wt{C}_n[\phi_{ij}]$ from
$\wt{C}_n[\R^m]$ to itself, where $\phi_{ij}$ are the coordinate
transformations defining $TM$.
\item Let $IC_e[M]$ be defined as in \refD{DT} but with 
$IC_n[M]$ replacing $IC_n(M)$.
\item Let $D_T[M]$ be defined
through the pull-back
$$
\begin{CD}
D_T[M] @>>> (IC_e[M])^{E_0} \\
@VVV  @VVV \\
C_{\#v_0}[M] @>>> (M)^{E_0}.
\end{CD}
$$
\end{enumerate}
\end{definition}

\begin{theorem}\label{T:closedstrata}
$C_T[M]$ is diffeomorphic to $D_T[M]$.
\end{theorem}

\begin{proof}
Though by definition $C_T[M]$ is the closure of $C_T(M)$ in $C_n[M]$,
it is also the closure of $C_T(M)$ in any closed subspace of $A_n[M]$,
and we choose to consider it as a subspace of $A_T[M]$.
The inclusion of $C_T(M)$ in $A_T[M]$ is compatible with fiber bundle
structures of these spaces over $(M)^{E_0}$.  For a general fiber bundle
$F' \to E' \to B'$ subspaces respectively of $F \to E \to B$, the 
closure $cl_{E}(E')$ may be defined by
first extending $E'$ to a bundle over $cl_{B}(B')$ (which may be done locally)
and then taking
the closures fiber-wise.  Our result follows from
this general statement, the definition 
of $C_{\#v_0}[M]$ as the closure of $C_{\#v_0}(M)$ in $A_{\#v_0}[M]$,
and the independence of the closure of the fibers $\wt{C}_i(\R^m)$
of $IC_i(M)$ in any $\wt{A}_i[\R^m]$.
\end{proof}

\subsection{Configurations in the line and associahedra} 

The compactification of configurations of points in the line is a fundamental case
of this construction.  The configuration spaces
$C_n(\R)$ and $C_n(\I)$ are disconnected, having one component for
each ordering of $n$ points. 
These different components each map to a
different component of $A_{n}[\R]$, because whether $x_i < x_j$ or
$x_i > x_j$ will determine a $+$ or $-$ for $u_{ij} \in S^0$. Let $C_n^o[\R]$
and $C_n^o[\I]$ denote the closure of the single component $x_1 < \cdots < x_n$.

The main result of this subsection is that $\wt{C}_n[\R]$ is Stasheff's associahedron $A_{n-2}$,
of which there is a pleasing description
of $A_n$ which we learned from Devadoss.  The truncation of a polyhedron
at some face (of any codimension) is the polyhedral subspace of points which
are of a distance greater than some sufficiently small epsilon from that 
face.  
We may define $A_n$ as a truncation of $\Delta^n$.  In the standard way, 
label the codimension one faces of $\Delta^n$ with elements of $\und{n+1}$.  
Call $S \subset \und{n+1}$ consecutive if
$i, j \in S$ and $i < k < j$ implies $k \in S$, and call a face of $\Delta^n$
consecutive if the labels of codimension one faces containing it are consecutive.
To obtain $A_n$, truncate the consecutive faces of $\Delta^n$, starting with
the vertices, then the edges, and so forth.

\skiphome{
\begin{center}
\begin{minipage}{10cm}
\begin{mydiagram}\label{F:assoc}\begin{center}
The third associahedron.
$$\includegraphics[width=6cm]{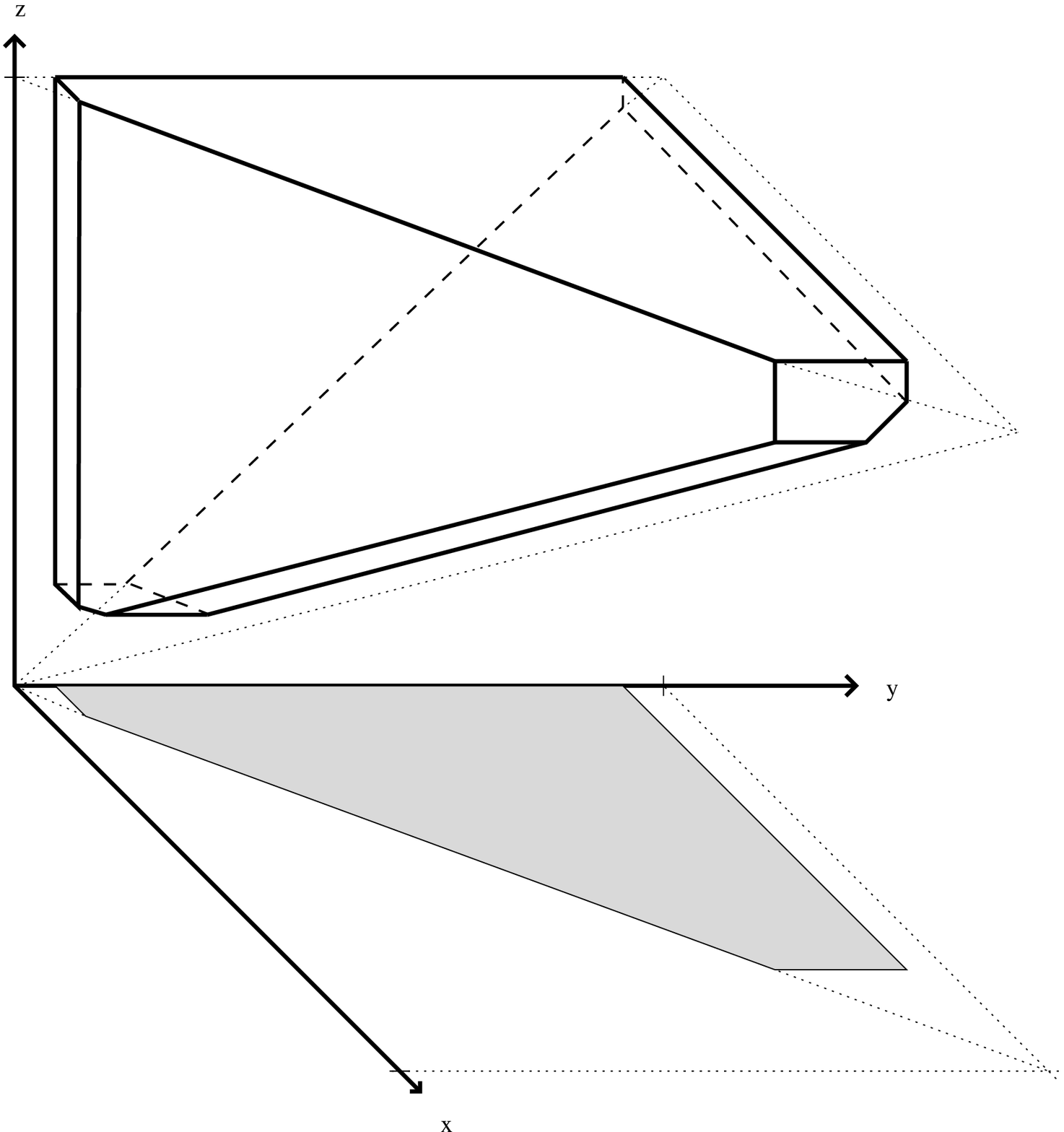} $$

\end{center}
\end{mydiagram}
\end{minipage}
\end{center}
}

We will use a more conventional definition of the associahedron below.  
Closely related to the associahedron is the following sub-category
of $\Psi_n$, whose minimal objects correspond to ways in which one can
associate a product of $n$ factors in a given order.

\begin{definition}
Let $\Psi^o_n$ denote the full sub-category of $\Psi_n$ whose objects
are $f$-trees such that the set of leaves over any vertex is consecutive
and such that the root vertex has valence greater than one.
\end{definition}

Note that any element of $\Psi^o_n$ has an embedding in the upper half plane with 
the root at $0$, in which the leaves occur in  order and which is unique up to isotopy.  
We may then drop the labels 
from such an embedding.

In applications to knot theory, we consider manifolds
with boundary which have two distinguished points in its
boundary, the interval $\I$ being a fundamental case.  

\begin{definition}
Given a
manifold $M$ with $y_0$ and $y_1$ in $\partial M$, 
let $C_n[M, \partial]$ be the closure in
$C_{n+2}[M]$ of the subspace of $C_{n+2}(M)$ of points of the form 
$(y_0, x_1, \ldots, x_n, y_1) \in (M)^{n+2}$.  
\end{definition}

\begin{theorem}\label{T:stashdiff}
Stasheff's associahedron $A_n$, $\wt{C}_{n+2}^o[\R]$ and $C_n^o[\I, \partial]$
are all diffeomorphic as manifolds with corners.
Moreover, their barycentric subdivisions are diffeomorphic to 
the realization (or order complex) of the poset $\Psi^o_n$.
\end{theorem}

\begin{proof}
It is simple to check that $\wt{C}_{n+2}^o[\R]$ and $C_n^o[\I, \partial]$
are diffeomorphic using \refL{sim}, and the fact that up to translation
and scaling, any $x_0 < x_1 < \cdots < x_{n+1} \subset \R$ has 
$x_0 = 0$ and $x_{n+1} = 1$.

Next, we analyze $\wt{C}_{n+2}^o[\R]$ inductively using \refC{wtCT} and
\refT{closedstrata}.  First note that because the $x_i$ are ordered and $x_0$ 
can never equal $x_{n+1}$, the category of strata of $\wt{C}_{n+2}^o[\R]$ is
$\Psi^o_{n+2}$.  For $n+2=3$,
$\wt{C}_{3}^o[\R]$ is a one-manifold whose interior is the open interval
$\wt{C}_3^o(\R)$ and which according to $\Psi^o_3$ has two distinct boundary
points, and thus must be an interval. For $n+2 = 4$, the stratification
according to
$\Psi^o_4$ and
\refT{closedstrata} dictate that there are five codimension-one boundary strata
each isomorphic to
$\wt{C}_3^o[\R]$, which we know inductively to be $\I$, and five vertices, each
being the boundary of exactly two faces, attached smoothly (with corners) to an
open two-disk, making a pentagon.

In general, $\wt{C}_{n+2}^o[\R]$ has an open $n$-ball for an interior and
faces $\wt{C}_T^o[\R]$
which inductively we identify as $(A_{\#v - 2})^{v \in V(T)}$ glued according
to the poset structure of $\Psi^o_{n+2}$ to make a boundary sphere, coinciding
with a standard definition of $A_n$ using trees \cite{Stas63}.

The last statement of the theorem follows from the general fact that if $P$ is a
polytope each of whose faces (including itself) is homeomorphic to a disk, then the
realization of the category of strata of $P$ is diffeomorphic to its barycentric
subdivision.
\end{proof}

In further work \cite{McSS03} we plan to show that the spaces $\wt{C}_n[\R^m]$
form an operad.  This construction unifies the associahedra and little disks operads, 
and was first noticed in \cite{GeJo95} and carried out in \cite{Mark99}.

To review some of the salient features of the structure of $C_n[M]$ in general,
it is helpful to think explictly about coordinates on $C_2^o[I, \partial]$.  On
its interior, suitable coordinates are $0 < x < y < 1$.  Three of the faces are
standard, corresponding to those for $\Delta^2$.  They are naturally labelled
$x = 0$, $y=1$ and $x=y$, and for example we may use $y$ as a coordinate on the
$x=0$ face, extending the coordinates on the interior.  The final two faces
are naturally labelled $0=x=y$ and $x=y=1$.   Coordinates on these faces which
extend interior coordinates would be $\frac{x}{y}$ and $\frac{1-y}{1-x}$,
respectively.

\skiphome{
\begin{center}
\begin{minipage}{10cm}
\begin{mydiagram}\label{F:assoc2}\begin{center}
The second associahedron, labelled by $\Psi^0_4$, with
labellings by associativity and coordinates also indicated.
$$\includegraphics[width=6cm]{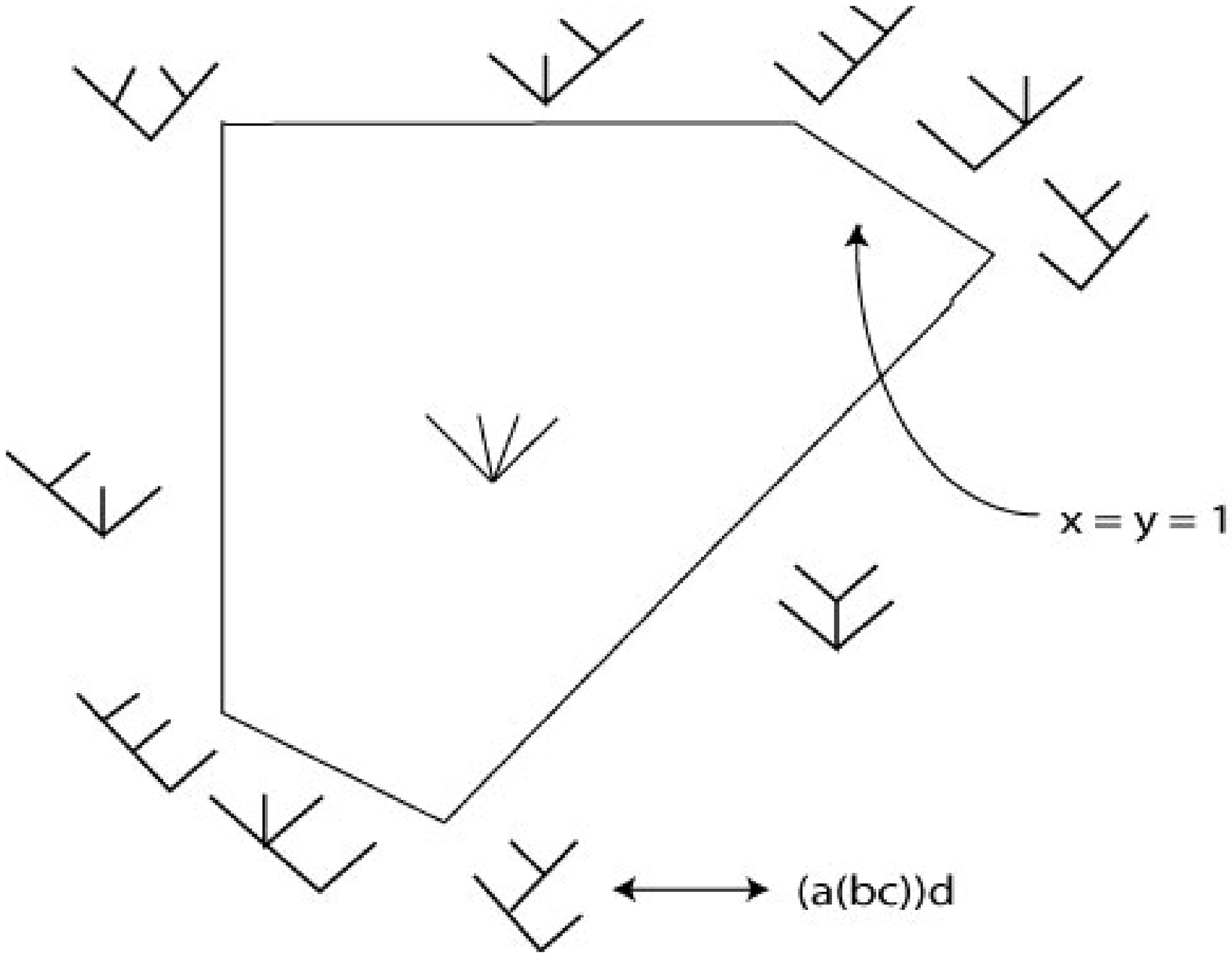} $$

\end{center}
\end{mydiagram}
\end{minipage}
\end{center}

}

\section{The simplicial compactification}\label{S:simplicial}

To simplify arguments, we assume $M$ is compact throughout this section.
Recall \refD{maindef} of $C_n \la M \ra$, which we call the simplicial
compactification.  For $M = \I$, we see that $C_n^o \la I \ra$ is 
the closure of $C_n^o(\I)$ in $\I^{\und{n}}$, which is simply $\Delta^n$.  
For general manifolds, we will see that $C_n \la M \ra$ is in some sense 
more complicated than $C_n[M]$.

Because the projection $P_A : A_{n}[M] \to A_n\la M \ra$  commutes 
with the inclusions of $C_n(M)$, 
\refL{projclose} says that $P_{A}$ sends $C_n[M]$ onto $C_n \la M \ra$
when $M$ is compact.

\begin{definition}
Let $Q_n \colon C_n[M] \to C_n \la M \ra$ be the restriction of $P_{A}$.
\end{definition}

The main aim of this section is to understand $Q_n$ and show that
it is a homotopy equivalence.  
From the analysis of \refL{lem1} we know that 
$(\pi_{ij})  : \wt{C}_n(\R^m) \to (S^{m-1})^{n(n-1)}$ is not injective for
configurations in which all points lie on a line. These collinear
configurations account for  all of the differences between
$C_n[M]$ and $C_n \la M \ra$.

\begin{lemma}\label{L:almostinj}
The map $Q_n$ is one-to-one except at points with some $x_{i_1} =
\cdots = x_{i_m}$ and $u_{i_h i_j} = \pm u_{i_k i_\ell}$ for any $h, j, k,
\ell$.  The preimages  of such points are diffeomorphic to a product of
$A_{m-2}$'s.
\end{lemma}

\begin{proof}
Conditions~\ref{d1} and \ref{d2} of \refT{cn[sub} say that in cases except for these,
the coordinates $d_{ijk}$ in $C_n[M]$ will be determined by the $x_i$ or $u_{ij}$  
coordinates.  In these cases, the $d_{ijk}$ are restricted in precisely
the same manner as for the definition of $\wt{C}_m^o(\R)$, which is
diffeomorphic to
$A_{m-2}$ by \refT{stashdiff}.
\end{proof}

Thus, the preimage of any point under $Q_n$ will be contractible, pointing
to the fact that $Q_n$ is a homotopy equivalence.  A small difficulty is that 
under $Q_n$,  points in the boundary of
$C_T[M]$ will be identified with points in its interior.  Moreover, there are
identifications made which lie only in the boundary of $C_T[M]$. 
We will first treat configurations in $\R^m$ up to the action of $Sim_k$, the building
blocks for the strata of $C_n[M]$.

\begin{definition}
Let $\wt{C}_n \la \R^m \ra$ be the closure of  $(\wt{\pi}_{ij})\left(\wt{C}_n(\R^m)\right)$
in  $\wt{A}_n \la \R^m\ra = (S^{m-1})^{C_2(\und{n})}$.
\end{definition}

The analogue of \refL{sim} does not hold in this setting, since as noted before 
$(\wt{\pi}_{ij})$ is not injective for collinear configurations.
We will see that $\wt{Q}_n : \wt{C}_n[\R^m] \to \wt{C}_n \la \R^m\ra$
is a homotopy equivalence by exhibiting $\wt{Q}_n$ as a pushout by an equivalence.  We
first state some generalities about fat wedges and pushouts.

\begin{definition}
Let $\{A_i \subseteq X_i\}$ be a collection of subpsaces indexed by $i$  in 
some finite $\mathcal{I}$.
Define the fat wedge of $\{X_i\}$ at $\{A_i\}$, denoted $\boxtimes^{\mathcal{I}}_{A_i} X_i$
or just $\boxtimes_{A_i} X_i$,
to be the subspace of $(x_i) \in (X_i)^{\mathcal{I}}$
with at least one $x_i$ in $A_i$.
\end{definition}

Suppose for each $i$ we have a map $q_i : A_i \to
B_i$ and let $Y_i$ be defined by the following pushout square:
$$
\begin{CD}
A_i @>>> X_i\\
@Vq_iVV  @V\bar{q_i}VV \\
B_i @>>> Y_i.
\end{CD}
$$

There is a map which we call $\boxtimes q_i$ from $\boxtimes_{A_i} X_i$ to
$\boxtimes_{B_i} Y_i$.  

\begin{lemma}\label{L:push}
With notation as above, 
if each $A_i \hookrightarrow X_i$ is a (Hurewicz) cofibration and each $q_i$ is a
homotopy equivalence, then $\boxtimes q_i$ is a 
homotopy equivalence.
\end{lemma}

\begin{proof}
First note that in a left proper model category,  if you have a diagram
$$
\begin{CD}
   B  @<<<  A  @>>>  C \\
    @VVV @VVV @VVV \\
     X @<<< Y @>>> Z,
\end{CD}
$$
 where the vertical maps are equivalences and at least one map
 on each of the horizontal levels is a cofibration, then the
 induced map of pushouts is an equivalence 
(see Theorem 13.5.4 in \cite{Hirs03}).   
The Hurewicz model category is left proper because every space is cofibrant
(see Theorem 13.1.3 in \cite{Hirs03}).

We prove this lemma by induction.  Let $\mathcal{I} = \{1, \ldots, n\}$.
Inductively define the diagram $D_j$ as
$$
\begin{CD}
     P_{j-1} \times X_j @<<<  P_{j-1} \times A_j  @>>>  (X_i)^{i<j} \times A_j \\
    @VVV @VVV @VVV \\
     Q_{j-1} \times Y_j @<<< Q_{j-1} \times B_j @>>> (Y_i)^{i<j} \times B_j,
\end{CD}
$$
where $P_{j-1}$ is the pushout of the top row of $D_{j-1}$ and 
$Q_{j-1}$ is the pushout of its bottom row.  Thus, $P_n = \boxtimes_{A_i} X_i$ and
$Q_n = \boxtimes_{B_i} Y_i$. 
The vertical maps of $D_j$ are homotopy equivalences by induction,
which is immediate if $n=1$.  The left most horizontal maps are cofibrations 
because the product of a Hurewicz cofibration with the identity map is
a Hurewicz cofibration (see Corollary 1 in \cite{Lill73}).  We apply the pushout
result above to get that $P_n \to Q_n$ is an equivalence.

\end{proof}

For analysis of $\wt{Q}_n$, recall that for $k=1$,
$\wt{\pi}_{ij}$ sends $\wt{C}_n[\R]$ to $S^0 = \pm 1$.

\begin{definition}
\begin{enumerate}
\item Let $\lambda_m$ be
the image of the map from $(S^{m-1}) \times \wt{C}_m[\R] \to
\wt{C}_m[\R^m]$ which sends $v, x$ to $(\pi_{ij}(x) \cdot v) \times (s_{ijk}(x))
\in \wt{A}_m[\R^m]$.  We call $\lambda_m$ the subspace of collinear points in
$\wt{C}_m[\R^m]$.

\item  Let $q_m$ denote the projection of $(S^{m-1}) \times \wt{C}_m[\R]$ onto
$S^{m-1} \times \Sigma_m$ by sending each component of $\wt{C}_m[\R]$ to a
point.  Let $R_m$ be defined as the pushout
\begin{equation}\label{Eq:lambda}
\begin{CD}
(S^{m-1}) \times \wt{C}_m[\R]  @>\lambda_m>> \wt{C}_m[\R^m]\\
@Vq_nVV   @VVV \\
S^{m-1} \times \Sigma_m      @>>>  R_m.
\end{CD}
\end{equation}

\end{enumerate}
\end{definition}

Note that $R_m$ maps to $\wt{C}_m \la \R^m\ra$, factoring $\wt{Q}_n$.  
We will see that this map is a homeomorphism on the image of $\wt{C}_m(\R^m)$ in
$R_m$, but not on its boundary strata.

\begin{definition}
\begin{enumerate}
\item By the analogue of \refT{closedstrata}, $\wt{C}_T[\R^m]$
is diffeomorphic to $\left(\wt{C_{\#v}}[\R^m]\right)^{V(T)}$.  Let $\lambda_T
\subset C_T[\R^m]$ be the fat wedge
$\boxtimes_{\lambda_{\#v}} \wt{C_{\#v}}[\R^m]$.

\item Let $L_T$ denote the fat wedge
$\boxtimes_{S^{m-1} \times \Sigma_{\#v}} R_{\#v}$ and let $q_T = 
\boxtimes q_{\#v} : \lambda_T \to L_T$.

\item Let $\bigcup_T \lambda_T$ denote the union of the $\lambda_T$ in
$\wt{C}_n[\R^m]$.   Let $\bigcup_T L_T$ denote the union of $L_T$ with
identifications $q_T(x) \sim q_{T'}(y)$ if $x \in \lambda_T$ is equal to
$y \in \lambda_{T'}$.  Let
$\bigcup_T q_T : \bigcup_T \lambda_T \to \bigcup_T L_T$ denote the projection
defined compatibly by the $q_T$.
\end{enumerate}
\end{definition}

\begin{theorem}\label{T:Proj}
The projection map $\wt{Q}_n  : \wt{C}_n[\R^m] \to \wt{C}_n \la \R^m\ra$ 
sits in a
pushout square
\begin{equation}\label{Eq:2}
\begin{CD}
\bigcup_T \lambda_T @>>> \wt{C}_n[\R^m] \\
@V\bigcup_T q_TVV  @VQ_nVV\\
\bigcup_T L_T @>>> \wt{C}_n \la \R^m\ra.
\end{CD}
\end{equation}
\end{theorem}

Before proving this theorem we deduce from it
one of the main results of this section.

\begin{corollary}
$\wt{Q}_n$ is a homotopy equivalence.
\end{corollary}

\begin{proof}
If we apply \refL{push} to the pushout squares of Equation~\ref{Eq:lambda}
which define the $R_{\#v}$, we deduce that $q_T$ is a homotopy equivalence. 
Because the identifications in $\bigcup_T L_T$ are essentially defined
through $\bigcup_T q_T$, we deduce that $\bigcup_T q_T$ is a homotopy
equivalence.  Because the inclusion $\bigcup_T \lambda_T \to \wt{C}_n[\R^m]$ is a
cofibration, we see that $\wt{Q}_n$ is a pushout of a homotopy equivalence through
a cofibration, and thus is a homotopy equivalence itself.
\end{proof}

\begin{proof}[Proof of \refT{Proj}]
Let $X$ denote the pushout of the first three spaces in the square of
Equation~\ref{Eq:2}.  First note that the composite 
$\wt{Q}_n \circ (\bigcup_T q_T)^{-1} : \bigcup_T L_T \to \wt{C}_n \la \R^m\ra$
is well defined, since choices of $(\bigcup_T q_T)^{-1}$ only differ in their
$d_{ijk}$ coordinates.  By the definition of pushout, $X$ maps to $\wt{C}_n \la
\R^m\ra$ compatibly with $\wt{Q}_n$.  We show that this map $F$ is a
homeomorphism. 

First, $F$ is onto because $\wt{Q}_n$ is onto.  The key is that by construction $F$
is one-to-one.   Away from $\bigcup_T \lambda_T$,
$\wt{Q}_n$ is one-to-one essentially by \refL{almostinj}.  The projection
$\wt{Q}_n$ is not one-to-one only on $x \in \wt{C}_m [\R^m]$ with some collections
of
$\{i_j\}$ such that $u_{i_h i_j} = \pm u_{i_\ell i_m}$.  But such an $x$ is in
$\lambda_{T(x)}$. The map $\wt{Q}_n \circ (\bigcup_T q_T)^{-1}$ is one-to-one since
distinct points in $\bigcup_T L_T$ will have distinct $u_{ij}$ coordinates when
lifted to $\bigcup_T \lambda_T$ which remain distinct in $\wt{C}_n \la \R^m\ra$.

Finally since it is a pushout of compact spaces, $X$ is compact.  All spaces in
question are subspaces of metric spaces.  Thus, since $F$ is a one-to-one map
between metrizable spaces whose domain is compact, it is a homeomorphism onto its
image, which is all of $\wt{C}_n\la \R^m\ra$.
\end{proof}

\begin{theorem}
The map $Q_n \colon C_n[M] \to C_n \la M \ra$ is a homotopy equivalence.
\end{theorem}

\begin{proof}
On the interior $C_n(M)$, $Q_n$ is a homeomorphism.

The effect of $Q_n$ on $C_T[M]$
for non-trivial $T$ is through restriction to $P_{\#v_0}$ on the base
$C_{\#v_0}[M]$.  Working fiberwise, we see that $Q_n$ takes each
fiber bundle  $\wt{C}_i[\R^m] \to IC_i[M] \to M$ and pushes out fiberwise to
get $\wt{C}_i \la \R^m \ra \to IC_i \la M \ra \to M$.
As $\#v_0 < n$, by induction and \refT{Proj},
$Q_n$ restricted to any $C_T[M]$ is a homotopy equivalence.  Since the inclusions
of $C_T[M]$ in each other are cofibrations, we can build a homotopy inverse
inductively and deduce that $Q_n$ is a homotopy equivalence.
\end{proof}

We also identify 
$C_n \la M \ra$ as a subspace of $A_n\la M \ra$ for purposes
of defining maps.  One approach to this identification would be to 
use surjectivity of $Q_n$ and \refT{cn[sub}, but there are relationships
between the $u_{ij}$ coordinates of points in $C_n[M]$
implied by their relationships in turn
with the $d_{ijk}$ coordinates.  It is simpler to do this coordinate
analysis of $C_n \la M \ra$ more directly, starting with $\wt{C}_n \la \R^m \ra$.

\begin{definition}
\begin{enumerate}

\item A point $(u_{ij}) \in \wt{A}_n \la \R^m \ra$ is anti-symmetric if $u_{ij} = -u_{ji}$.

\item A circuit, or $k$-circuit, in $S$ is a collection 
$\{ i_1 i_2,$ $i_2 i_3$, \ldots, $i_{k-1} i_k \}$
of elements of $\binom{S}{2}$ for some indexing set $S$.  
Such indices label a path in the complete
graph on $S$.  A circuit is a loop, or $k$-loop,  if $i_k = i_1$.
A circuit is straight if  it does not contain any loops.
The reversal of a circuit is the circuit $i_k i_{k-1}, \ldots, i_2 i_1$.

\item  A point $(u_{ij}) \in \wt{A}_n \la \R^m \ra$ is three-dependent if for $3$-loop
 $L$  in $\n$ we have $\{u_{ij}\}_{ij \in L}$ is non-negatively dependent.
 
\item  If $S$ has four elements and is ordered 
we may associate to a straight $3$-circuit $C = \{ ij, jk, k\ell \}$
a permutation of $S$ denoted $\sigma(C)$ which orders $(i, j, k, \ell)$.
A complementary $3$-circuit $C^*$ is a circuit, unique up
to reversal, which is comprised of the three pairs of indices not in $C$.

\item A point in $\wt{A}_n \la \R^m \ra$ is four-consistent if for any $S \subset \n$ of 
cardinality four and any $v, w \in S^{m-1}$
we have that 
\begin{equation} \label{4consist}
\sum_{C \in \mathcal{C}^3(S)} (-1)^{|\sigma(C)|}
\left( \prod_{ij \in C} u_{ij} \cdot v \right)\left(\prod_{ij \in C^*} u_{ij} \cdot w \right) = 0,
\end{equation}
where $\mathcal{C}^3{S}$ is the set of straight $3$-circuits modulo reversal and
${|\sigma(C)|}$ is the sign of $\sigma(C)$.
\end{enumerate} 
\end{definition}

One may view anti-symmetry as a dependence condition for two-loops of indices.

\begin{lemma}\label{L:3d4c}
The image of $(\wt{\pi}_{ij})$ is anti-symmetric, three-dependent and four-consistent.
\end{lemma}

\begin{proof}
Let $(x_i)$ be a coset representative in $\wt{C}_n(\R^m)$ and $u_{ij} = \pi_{ij}\left((x_i)\right)$.
Three-dependence is immediate as $(x_i - x_j) + (x_j - x_k) + (x_k - x_i) = 0$,
so $||x_i - x_j|| u_{ij} + ||x_j - x_k|| u_{jk} + ||x_k - x_i|| u_{ki} = 0$.

Four-consistency is more involved.  We start in the plane, and we work projectively letting
$a$ be the slope of $u_{12}$, $b$ of $u_{34}$, $c$ of $u_{13}$, $d$ of $u_{24}$,
$e$ of $u_{14}$ and $f$ of $u_23$.  If $a = \infty$ and $b = 0$, so up to translation we arrange for  $x_1$ and $x_2$
on the $y$-axis and $x_3$ and $x_4$ on the $x$-axis, we observe that $cd = ef$.  
To lift our assumptions on $a$ and $b$, we use the fact that any linear fractional
transformation of slopes is induced by a linear transformation of the $(x_i)$, and thus 
preserves the slopes which come from some $(x_i)$.  
We apply the transformation $t \mapsto \frac{t-b}{t-a}$, which sends
 $a$ to $\infty$, $b$ to $0$, and the equation $cd - ef = 0$ to the equation 
$$(c-b)(d-b)(e-a)(f-a) - (c-a)(d-a)(e-b)(f-b) = 0.$$
The resulting quartic is divisible by $(a-b)$.  Carrying out the division we get the more
symmetric cubic
$$ab(-c -d + e + g) + (a+b)(cd-ef) + [ef(c+d) - cd(e+f)] = 0.$$
Now recalling that for example $a = \frac{u_{12} \cdot e_2}{u_{12} \cdot e_1}$ where
$\{ e_1, e_2 \}$ is the standard basis of the plane, we clear
denominators and find that we have precisely Equation~\ref{4consist}, in the case of the
plane where $v = e_2$ and $w = e_1$.
To deduce the case of general $v$ and $w$ in the plane, we 
simply change to the $v,w$ basis when $v$ and $w$ are independent.
The case of dependent $v$ and $w$ follows by continuity.
Finally, we invoke the fact that
 the dot product of $u_{ij}$ with $v$ is the same as that
of the projection of $u_{ij}$ onto the plane spanned by $v$ and $w$ to deduce
the general case from the planar case.
\end{proof}

\begin{lemma}\label{L:5gives6}
If $(u_{ij}) \in \wt{A}_4 \la \R^m \ra$ is four-consistent, then any five of the
$u_{ij}$ determines the sixth.
\end{lemma}

\begin{proof}
The four-consistency condition, Equation~\ref{4consist}, is multilinear in each variable,
and the terms are all nonzero for $v,w$ in the complement of the hyperplanes orthogonal
to the $\{u_{ij}\}$.  Thus for generic $v$ and $w$ one is
able to determine the ratio $\frac{u_{k \ell} \cdot v}{u_{k \ell} \cdot w}$ from knowing all other 
$u_{ij}$.  A unit vector is determined by such ratios, in fact 
needing only the ratios between pairs of vectors in some basis.
\end{proof}

\begin{remark}
If  $u_{ik}$, $u_{i \ell}$ and $u_{jk}$ are independent, then the four-consistency
condition follows from three-dependency, as $u_{k \ell}$ 
must be the intersection of the plane through the origin,
$u_{ik}$ and $u_{i \ell}$, and the plane
through the origin, $u_{jk}$ and $u_{j\ell}$.
\end{remark}

\begin{theorem}\label{T:cnlasub}
$\wt{C}_n\la\R^m\ra$ is the subspace of anti-symmetric, three-dependent, four-consistent points
in $\wt{A}_n \la \R^m \ra$.
\end{theorem}

\begin{proof}
Let $DC$ be the subspace of $\wt{A}_n \la \R^m \ra$ of anti-symmetric three-dependent,
four-consistent points. We proved in \refL{3d4c} points in the image of $(\pi_{ij})$ are
in $DC$.  Moreover, $DC$ is closed, since anti-symmetry, dependence of vectors and four-consistency
are closed conditions.  We establish the theorem by constructing points in the
image of $(\pi_{ij})$ arbitrarily close to any point in $DC$, in a manner reminiscent of
the maps $\nu^T$.

To a point $u = (u_{ij}) \in \wt{C}_n\la\R^m\ra$ we associate an exclusion relation and thus
using \refD{extree} an $f$-tree $T(u) \in \wt{\Psi}_{\n}$, by saying that
$i$ and $j$ exclude $k$ if $u_{ik} = u_{jk} \neq \pm u_{ij}$.  It is immediate that this satisfies
the first axiom of an exclusion relation.  To check the second axiom, namely transitivity, 
we also assume that $j$
and $k$ exclude $\ell$ so that $u_{j\ell} = u_{k \ell} \neq \pm u_{jk}$.  We use four-consistency
with $v$ chosen to be orthogonal to $u_{ik}$ but not $u_{ij}$ or $u_{j \ell}$
(we may assume that we are
not in $\R^1$, for in that case the exclusion relation would be empty automatically) and
$w$ orthogonal to $u_{j \ell}$ but not $u_{jk}$.  
All but one of the twelve terms of 
Equation~\ref{4consist} are automatically zero, with the remaining term being
$$ (u_{ij} \cdot v) (u_{j \ell} \cdot v)(u_{\ell k} \cdot v)
                      (u_{jk} \cdot w)(u_{ki} \cdot w)(u_{i \ell} \cdot w ).$$
All of these factors are non zero by construction except for $u_{i \ell} \cdot w$.  We deduce
that $u_{i \ell}$ is orthogonal to all vectors $w$ which are orthogonal to  $u_{j \ell}$, which we
recall is not orthogonal to $u_{ij}$ so that 
$u_{i \ell} = \pm u_{j \ell} \neq \pm u_{ij}$.  By the non-negativity of coefficients in
our three-dependence condition we in 
fact have $u_{i \ell} = u_{j \ell}$ so that $i$ and $j$ exclude $\ell$, as needed.

We next show that if there are no exclusions regarding indices, we can construct 
a non-continuous inverse $\rho_n$ to $(\pi_{ij})$.  
If there are no exclusions, then either:
\begin{enumerate} 
\item for any subset of indices $S$
there is a $k$ such that $u_{ik} \neq u_{jk}$ for some $i,j \in S$, \label{one}
\item or, inductively, all of the $u_{ij}$ are equal up to a sign. \label{two}
\end{enumerate}
In case~(\ref{two}) let $v = u_{12}$ and
define a total ordering on $\n$ by $i < j$ if $u_{ij} = v$.   
Define $(x_i) = \rho_n(u_{ij})$  by setting 
$x_i = s_i v$, where $s_i > s_j$ when $i>j$, $\Sigma s_i = 0$ and 
max$\{ |s_i| \} = 1$, noting that there is not a unique choice of $s_i$'s.  
In case~(\ref{one}), we  begin by setting  $x_1 = 0$ and $x_2 = u_{12}$. 
Once $\{x_i\}$ for $i \in S$ has been determined, 
choose $k$ as in (\ref{one}).  Let $R_{ik}$ be the positive ray 
through $x_i$ in the direction of $u_{ik}$.  Define $x_k$ as the
intersection of  of  $R_{ik}$ with $R_{jk}$, with $i$ and $j$ as in (\ref{one}).  
These rays intersect for they are coplanar by three-dependence.  
They are not parallel because $u_{ik} \neq  \pm u_{jk}$, and the intersection
of the two lines containing them lies in each ray
because the dependence of $u_{ij}$, $u_{jk}$ and $u_{ki}$ is through
positive coefficients.

To show that this definition of $x_k$ is independent of which 
$i$ and $j$ are chosen, we show that we may replace $j$ by some
$\ell$.  Letting $S = \{i,j,k, \ell\}$, we observe that both the image of
these $(x_r)_{r \in S}$ under $(\pi_{pq})_{p,q \in S}$ and our given 
$(u_{pq})_{p,q \in S}$ are four-consistent, by \refL{3d4c} and
by assumption respectively.  By construction, $\pi_{pq}((x_r)_{r \in S}) = u_{pq}$
in all cases except perhaps when $p = k$, $q = \ell$.  But \refL{5gives6}, 
implies this equality in this last case.
Thus  $x_k = R_{ik} \cap R_{jk}$ is on the ray $R_{\ell k}$, which in turn implies
that $x_k$ is also $R_{ik} \cap R_{\ell k}$.   

In this way we construct $x_i$ for all $i \in \n$.
By scaling and translating, we choose 
the coset representative for $\rho_n((u_{ij}))$ to be the $(x_i)$ whose 
center of mass is the origin and such that max$\{|| x_i ||\} = 1$.

Finally, we use these $\rho_i$ to construct a configuration in $C_n(\R^m)$ which maps
arbitrarily close to any given $u = (u_{ij}) \in DC$.  For every vertex $v \in T(u)$
we choose one leaf lying over each edge of $v$ and let $S_v$ be the set of
their indices.  There are no exclusions among the indices in $S_v$,
so we let $(x_e)_{e \in E(v)} = \rho_{\# S_v}((u_{ij})_{i,j \in S_v})$.  Given
an $\varepsilon < 1$ let $x_i(\varepsilon) = \sum_{e \in R} \varepsilon^{h(e)} x_e$,
where $R$ is the root path of the $i$th leaf, $e$ is an edge in that root path, and $h(e)$
is the number of edges in $R$ between $e$ and the root vertex.  By construction,
independent of all of the choices made, the image of $(x_i(\varepsilon))$
under $(\pi_{ij})$ approaches $(u_{ij})$ as $\varepsilon$ tends to zero.
\end{proof}

\begin{corollary}
$C_n \la M \ra$ is the subspace of $(x_i) \times (u_{ij}) \in A_n \la M \ra$ such that 
\begin{enumerate}
\item if $x_{i} \neq x_k$ then $u_{ij} = \frac{x_i - x_k}{||x_i - x_k||}$;
\item the $(u_{ij})$ are anti-symmetric, three dependent and four-consistent;
\item  If $x_i = x_j$, then $u_{ij}$ is tangent to $M$ at $x_i$.
\end{enumerate}
\end{corollary}

We may decompose $\wt{C}_n \la M \ra$ through our construction of associated trees $T(u)$.
These strata are manifolds,  but $C_n \la M \ra$ is not a manifold with corners.   
The singularity which arises is akin to that
which occurs when a diameter of a disk gets identified to a point.  We will not
pursue the matter further here.

%\begin{figure}\label{Qn}
%
%{
%\psfrag{A}{$Q_n$}
%\psfrag{B}{$\lambda_T$}
%\psfrag{C}{$L_T$}
%\psfrag{D}{$C_T\la M \ra$}
%\psfrag{E}{$C_T[M]$}
%$$\includegraphics[width=9cm]{quotient.eps}$$}
%\caption{Singularity caused by $Q_n$.}
%\end{figure}

\section{Diagonal and projection maps}  

As we have seen, the compactifications $C_n[M]$ and $C_n \la M \ra$
are functorial with respect to embeddings of $M$.
In this section we deal with projection and diagonal maps, 
leading to functorality with respect to $n$, viewed as the set $\und{n}$.

Our goal is to construct maps for $C_{\#S}[M]$ and $C_{\#S}\la M \ra$
which lift the canonical maps on $M^S$.  We start with the straightforward case of
projection maps.
If $\sigma : \und{m} \to \und{n}$ is an inclusion of sets, recall \refD{funct}
that $p^M_{\sigma}$ is the projection onto coordinates in the image of
$\sigma$.  

\begin{proposition}\label{P:projmaps}
Let $\sigma : \und{m} \to \und{n}$ be an inclusion of finite sets. There are
projections $C_\sigma$ from $C_n[M]$ onto $C_m[M]$ and from $C_n \la M
\ra$ onto $C_m \la M \ra$ which are consistent with each other, with 
$p_{\sigma}^M$, and its restriction to $C_n(M)$.
\end{proposition}

\begin{proof}
The inclusion $\sigma$ gives rise to maps from
$C_i(\sigma) : C_i(\und{m}) \to C_i(\und{n})$.  
We project $A_{n}[M]$ onto 
$A_{m}[M]$ through $P_\sigma = p_\sigma^M \times p_{C_2(\sigma)}^{S^{m-1}}
\times p_{C_3{\sigma}}^{\I}$.  

Because $P \circ \alpha_n = \alpha_m$ and all spaces in question are
compact, we apply  \refL{projclose} to see that $P_\sigma$ sends
$C_n[M]$ onto $C_m[M]$, extending the
projection from $C_n(M)$ to
$C_m(M)$. By construction, $P_\sigma$ commutes with $p_\sigma^M$,
which is its first factor.  

The projection for $C_n \la M \ra$ is entirely analogous, defined as
the restriction of the map $P'_\sigma = p_\sigma^M \times
p_{C_2(\sigma)}^{S^{m-1}} : A_n \la M \ra \to A_m \la M \ra$.
We leave the routine verification that $P'$ commutes with all maps
in the statement of the theorem to the reader.
\end{proof}

An inclusion $\sigma : \und{m} \to \und{n}$ gives rise to a functor 
$Ex_\sigma : Ex(\und{n}) \to Ex(\und{m})$ by throwing out any exclusions involving
indices not in the image of $\sigma$.  The corresponding ``pruning'' functor for
trees, $\Psi_\sigma : \Psi_{\und{n}} \to \Psi_{\und{m}},$  is defined by removing
leaf vertices and edges whose label is not in the image of $\sigma$, replacing any
non-root bivalent vertex along with its two edges with a single edge, and
removing any vertices and edges which have all of the leaves above them removed.  

\begin{proposition}
$C_\sigma$ sends $C_T[M]$ to $C_{\Psi_\sigma(T)}[M]$.
\end{proposition}

\begin{proof}
The effect of $C_\sigma$ is to omit indices not in the image of $\sigma$, so its
effect on exclusion relations is precisely $Ex_\sigma$.  There is a univalent root
vertex for the tree associated to $C_\sigma(C_T[M])$ if and only if all indices
$j$ for which $x_j \neq x_i$ have been omitted, which happens precisely when all
leaves in $T$ except for those over a single root edge have been pruned.
\end{proof}

If $\sigma : \und{m} \to \und{n}$ is not injective, it is more
problematic to construct a corresponding map $C_n [M] \to C_m[M]$.
Indeed, $p_\sigma : M^{\und{n}} \to M^{\und{m}}$ will not send $C_n(M)$
to $C_m(M)$, since the image of $p_\sigma$ will be some diagonal subspace
of $M^{\und{m}}$ and the diagonal subspaces are precisely what are
removed in defining $C_n(M)$.  One can attempt to define diagonal maps by
``doubling'' points, that is adding a point to a configuration which is
very close to one of the points in the configuration, but such
constructions are non-canonical and will never satisfy the identities which
diagonal maps and projections together usually do.  But, the
doubling idea carries through remarkably well for compactified
configuration spaces where one can ``double infinitesimally''.  From the viewpoint
of applications in algebraic topology, where projection and diagonal maps are
used frequently, the diagonal maps for compactifications of 
configuration spaces should be of great utility.

Reflecting on the idea of doubling a point in a configuration, we see that
doing so entails choosing a direction, or a unit tangent vector, at that
point.  Thus we first incorporate tangent vectors in our constructions.
Recall that we use $STM$ to denote the unit tangent bundle (that is, the
sphere bundle to the tangent bundle) of $M$.
% and let $p$ be the projection of $STM$ onto $M$.

\begin{definition}
If $X_n(M)$ is a space with a canonical map to $M^{\und{n}}$, define
$X'_n(M)$ as a pull-back as follows:
$$
\begin{CD}
X'_n(M)  @>>> (STM)^{\und{n}} \\
@VVV @VVV\\
X_n(M) @>>> M^{\und{n}}.
\end{CD}
$$
If $f_n : X_n(M) \to Y_n(M)$ is a map over $M^{\und{n}}$, let $f'_n 
:  X'_n(M) \to Y'_n(M)$ be the induced map on pull-backs.
\end{definition}

\begin{lemma}
$C_n'[M]$ is the closure of the image of $\alpha_n' : C'_n(M) \to
A_n'[M]$.  Similarly, $C_n' \la M \ra$ is the closure of the image of
$\beta_n'$.
\end{lemma}

\begin{proof}
We check that $cl_{A'_n[M]}\left(\alpha'_n(C'_n(M))\right)$ satisfies
the definition of $C'_n[M]$ as a pull-back by applying
\refL{projclose} with $\pi$ being the projection from 
$A'_n[M]$ to $A_n[M]$  and $A$ being the subspace $\alpha_n(C'_n(M))$. 
The proof for $C'_n \la M \ra$ proceeds similarly. 
\end{proof}

We may now treat both diagonal and projection maps for $C'_n \la M \ra$.
Starting with $M = \R^m$, note that 
$A'_n\la \R^m\ra = (\R^m\times S^{m-1})^{\und{n}} \times
(S^{m-1})^{C_2(\und{n})}$, which is canonically diffeomorphic to  $(\R^m)^{\und{n}}
\times (S^{m-1})^{\und{n}^2}$, as we let $u_{ii}$ be the unit tangent vector
associated to the $i$th factor of $\R^m$.  

\begin{definition}
Using the product decomposition above and considering $M$ as a 
submanifold of $\R^m$, define $A_\sigma :
A'_n\la \R^m\ra \to A'_m \la \R^m\ra$ as 
$p^{\R^m}_\sigma \times p^{S^{m-1}}_{\sigma^2}$
and let $F_\sigma$ be the restriction of $A_\sigma$ to $C_n\la M \ra$.
\end{definition}

\begin{proposition}\label{P:Fsigma}
Given  $\sigma : \und{m} \to \und{n}$ the induced map $F_\sigma$ sends 
$C'_n\la M \ra$ to $C'_m \la M \ra$ and commutes with $p_\sigma^{STM}$.
\end{proposition}

\begin{proof}
To see that the image of $F_\sigma$ lies in $C_m'\la M \ra$, it
suffices to perform the routine check 
that its projection to $A_m \la M \ra$ satisfies the conditions of
\refT{cnlasub} using the fact that the domain of $F_\sigma$, namely $C'_n \la M
\ra$, satisfies similar conditions.   Let $(x_i) \times (u_{ij})$ be
$F_\sigma\left((y_\ell) \times (v_{\ell m})\right)$ so that $x_i = y_{\sigma(i)}$
and
$u_{ij} = v_{\sigma(i) \sigma(j)}$.  

Looking at the first condition of \refT{cnlasub}, 
$x_i \neq x_j$ means $y_{\sigma(i)} \neq y_{\sigma(j)}$.  By
\refT{cnlasub} applied to $C_n \la M \ra$ we have that $v_{\sigma(i) \sigma(j)}$
is the unit vector from $y_{\sigma(i)}$ to $y_{\sigma(j)}$, which implies the
corresponding fact for $u_{ij}$.  Similarly, that the $(u_{ij})$
are anti-symmetric, three-dependent and four-consistent follows
mostly from the corresponding statements for the $(v_{\ell m})$.  Cases
which do not follow immediately in this way, such as three-dependence 
when the indices $i$, $j$ and
$k$ are not distinct, are degenerate and thus straightforward to verify, 
as for example in this case two of these vectors will be equal up to sign.
We leave such verification to the reader.
\end{proof}

Let $\mathcal{N}$ denote the full subcategory of the category of sets generated by 
the $\und{n}$.

\begin{corollary}
Sending $\und{n}$ to $C_n' \la M \ra$ and $\sigma$ to $F_\sigma$ defines
a contravariant functor from $\mathcal{N}$ to spaces.
\end{corollary}

\begin{proof}
We check that $F_{\sigma \circ \tau} = F_\sigma \circ F_\tau$.  This follows from
checking the analogous facts for $p_\sigma$ and $p_{\sigma^2}$, which are
immediate.
\end{proof}

Let $[n] = \{ 0, \cdots, n \}$, an ordered set given the standard ordering of integers.
Recall the category $\Delta$, which has one object for each nonnegative $n$ and
whose morphisms are the non-decreasing ordered set morphisms between the $[n]$.  A
functor from $\Delta$ to spaces is called  a cosimplicial space.  There is a canonical 
cosimplicial space often denoted $\Delta^\bullet$ whose $n$th object is $\Delta^n$.
To be definite we coordinatize $\Delta^n$ by
$0 = t_0 \leq t_1  \leq \cdots \leq t_n \leq t_{n+1} = 1$, and label its vertices by
elements of $[n]$ according to the number of $t_i$ equal to one.
The structure maps for this standard object
are the linear maps extending the maps of vertices as 
sets.  On coordinates, the linear map corresponding to some $\sigma: [n] \to [m]$ 
sends  $(t_i) \in \Delta^n$ to $(t_{\sigma^*(j)}) \in \Delta^m$ where $n-\sigma^*(j)$ 
is the  number of $i \in [n]$ such that $\sigma(i) < m-j$.

The following
corollary gives us another reason to refer to $C_n \la M \ra$ as the simplicial
compactification of $C_n(M)$.  For applications we are interested in a
manifold $M$ equipped with one inward-pointing tangent vector $v_0$ and one
outward-pointing unit tangent vector $v_1$ on its boundary. Let $C_n' \la M, \partial
\ra$ denote the subspace of $C_{n+2}' \la M, \partial \ra$ whose first projection onto
$STM$ is $v_0$ and whose $n+2$nd projection is $v_1$.   Let 
$\phi : \Delta \to \mathcal{N}$ be the functor which sends $[n]$ to $\und{n+1}$ 
and relabels the morphism accordingly.

\begin{corollary}
The functor which sends $[n]$ to $C'_n \la M, \partial  \ra$ and $\sigma: [n] \to [m]$ to
the restriction of $p_{\tau}$ to $C'_n \la M, \partial \ra$ where $\tau : [m+1] \to [n+1]$ 
is the composite $  \phi \circ \sigma^*  \circ \phi^{-1}$
defines a  cosimplicial space.  
\end{corollary}

This cosimplicial space models the space of knots in $M$ \cite{Sinh02}.

For $C_n'[M]$, projection maps still work as in
\refP{projmaps}, but diagonal maps are less canonical and more involved to
describe.  We restrict to a special class of diagonal maps for simplicity.

\begin{definition}
Let $\sigma_i : \und{n+k} \to \und{n}$ be defined by letting
$K_i = \{i, i+1, \ldots, i+k\}$ and setting 
$$
\sigma_i(j) = 
\begin{cases}
j & j<K_i \\
i & j \in K_i \\
j-k & j > K_i.
\end{cases}
$$
\end{definition}

We must take products with associahedra in order to account for all possible
diagonal maps.

\begin{definition}
\begin{enumerate}
\item Define $\iota_i : I^{C_3(\und{n})} \times A_{k-1} \to I^{C_3(\und{n+k})}$
by recalling that $A_{k-1} \cong \wt{C_{k+1}(\R)} \subset I^{C_3(\und{k})}$ and
sending $(d_{j\ell m})^{C_3(\und{n})} \times (e_{j\ell m})^{C_3(\und{k})}$ to
$(f_{j\ell m})^{C_3(\und{n+k})}$ with 
$$
f_{j\ell m} =
\begin{cases}
d_{\sigma_i(j,\ell,m)} & {\text{if at most one of}} \; j, \ell, m \in K_i \\
0 & {\text{if}} \; j,\ell \in K_i \; {\text{but}} \; m \notin K_i \\
1 & {\text{if}} \; \ell, m \in K_i \; {\text{but}} \; j \notin K_i \\ 
\infty & {\text{if}} \; j,m \in K_i \; {\text{but}} \; \ell \notin K_i \\
e_{j-i, \ell-i, m-i} & {\text{if}} \; j, \ell, m \in K_i.
\end{cases}
$$

\item Let $D_{i,k} : A_n'[M] \times A_{k-1} \to A_{n+k}' [M]$ be the product of
$A_{\sigma_i} : A_n' \la M \ra \to A_{n+k}' \la M \ra$ with $\iota_i$.
Let $\delta^i_k$ denote the restriction of $D_i$ to $C_n'[M] \times A_{k-1}$.
\end{enumerate}
\end{definition}

\begin{proposition}\label{P:deltai}
$\delta^i_k$ sends $C_n'[M] \times A_{k-1}$ to $C_{n+k}'[M] \subset A'_{n+k}[M]$.
\end{proposition}

As with \refP{Fsigma}, the proof is a straightforward checking that the image of
$\delta^i_k$ satisfies the conditions of \refT{cn[sub}.  One uses the fact that 
$C_n'[M]$ satisfies those conditions, along with the definition of $\iota_i$.  We
leave a closer analysis to the reader.  

By analysis of the exclusion relation, we see
that the image of $\delta^i_k$ lies in
$C'_S[M]$ where $S$ is the tree with $n+k$ leaves, where leaves with labels in
$K_i$ sit over the lone one internal vertex, which is initial for the $i$th root
edge.  In general, $\delta^i_k$ sends $C'_T[M]$  to
$C'_{T'}[M]$, where $T'$ is obtained from $T$ by adding $k+1$ leaves to $T$, each
of which has the $i$th leaf as its initial vertex.

We set $\delta^i = \delta^i_1 : C_n'[M] \to C_{n+1}'[M]$, and note that these act as
diagonal maps.  One can check that composing this with the projection down back
to $C_n'[M]$ is the identity.  Unfortunately, $\delta^i \delta^i \neq
\delta^{i+1} \delta^i$ -- see Figure~\ref{deltas} -- so that the $C_n'[M]$ do not form a
cosimplicial space.  But note that our $\delta^2$, when we restrict
$A_1$ to its boundary, restricts to these two maps and thus provides a
canonical homotopy between them.  In fact \refP{deltai} could be used to make an
$A_\infty$ cosimplicial space, but it is simpler to use the $C_n' \la M \ra$
if possible.

\skiphome{
\begin{center}
\begin{minipage}{10cm}
\begin{mydiagram}\label{deltas}\begin{center}
An illustration that $\delta^2 \delta^2 \neq \delta^3 \delta^2$.
\psfrag{D2}{$\delta^2$}
\psfrag{D22}{$\delta^2$}
\psfrag{D23}{$\delta^3$}
$$\includegraphics[width=9cm]{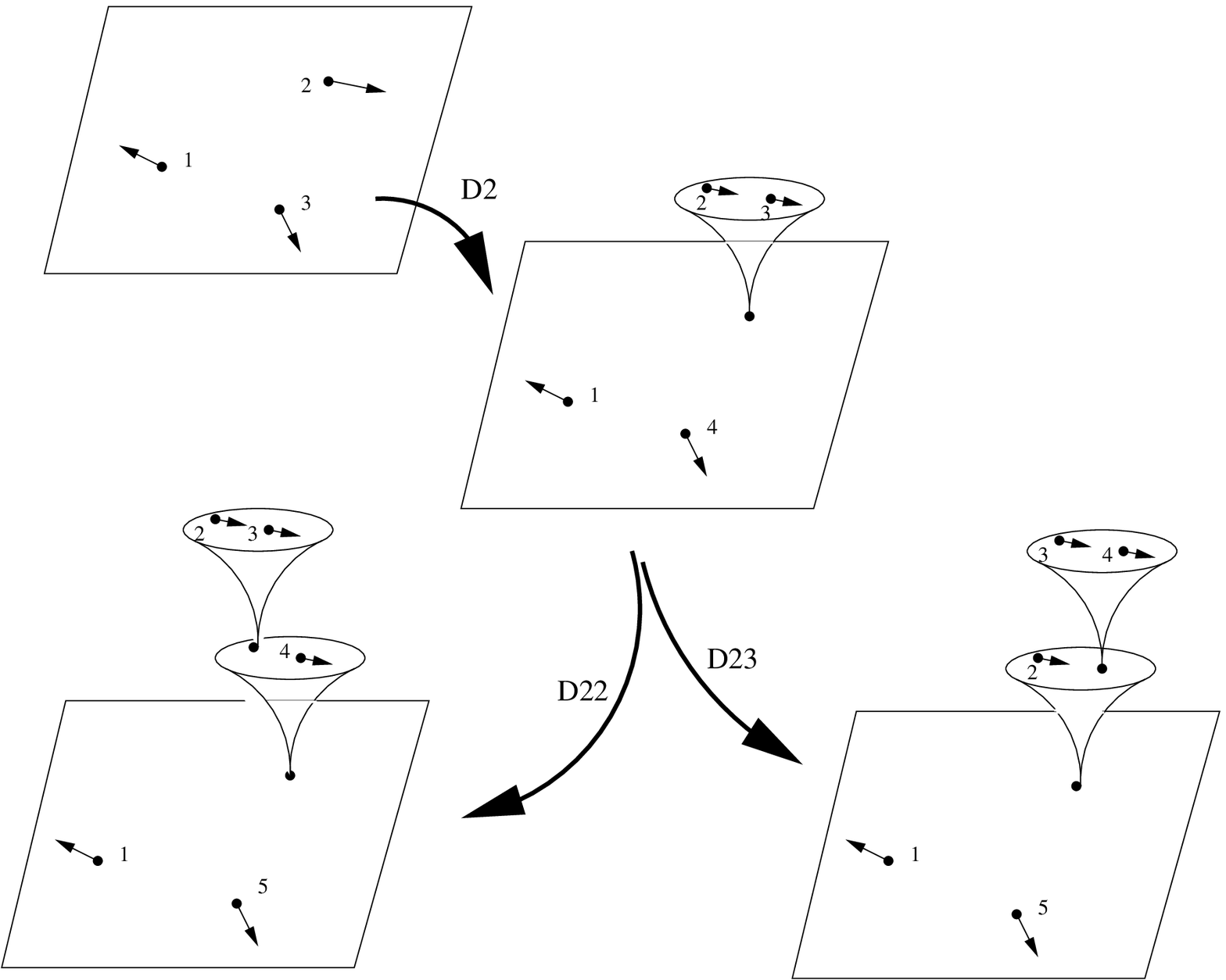}$$
\end{center}
\end{mydiagram}
\end{minipage}
\end{center}

}


\begin{thebibliography}{10} 

\bibitem{Axel94} S. Axelrod and I. Singer,   
Chern-Simons perturbation theory, II.  
{\em Jour. Diff. Geom.} {\bf 39} (1994), no. 1,
173--213.   

\bibitem{BGRT02} D/  Bar-Natan, S. Garoufalidis, L. Rozansky, and D. Thurston. The Aarhus 
integral of rational homology 3-spheres I. {\em  Selecta Math. (N.S.)}  {\bf 8}  
(2002),  no. 3, 315--339.

\bibitem{Bott94} R. Bott and C. Taubes,  \newblock On the self-linking
of  knots. Topology and physics. \newblock {\em J. Math. Phys.} {\bf
35} (1994), no. 10, 5247--5287. 

\bibitem{BCSS03} R. Budney, J. Conant, K. Scannell and D. Sinha,
New perspectives on self-linking.  To appear in {\em Advances in Mathematics}.

\bibitem{CCL02} A. Cattaneo, P. Cotta-Ramusino, and R Longoni, Configuration spaces and 
Vassiliev classes in any dimension.  {\em Algebr. Geom. Topol.}  {\bf 2}  (2002), 949--1000.

\bibitem{DCPr95} C. De Concini and C. Procesi,  Wonderful models of subspace
arrangements.  {\em Selecta Mathematica} (NS), {\bf 1}, (1995), 459--494.

\bibitem{Deva02} S. Devadoss, 
A space of cyclohedra. 
{\em Discrete Comput. Geom.} {\bf 29} (2003), no. 1, 61--75.

\bibitem{FH01} E. Faddell and S. Husseini,  {\em Geometry and topology
of configuration spaces}.  Springer, 2001. 

\bibitem{FeKo03} E. Feitchner and D. Kozlov,  Incidence combinatorics of
resolutions, preprint 2003.

\bibitem{Fult94} W. Fulton and R. MacPherson, 
Compactification  of configuration spaces.   {\em
Annals of Mathematics} {\bf  139} (1994), 183--225. 

\bibitem{Gaif03} G. Gaiffi,  Models for real subspace arrangements and stratified
manifolds.   {\em Int. Math. Res. Not.}  (2003)  no. 12, 627--656.

\bibitem{GeJo95} E. Getzler and J. Jones,  Operads, homotopy algebra and iterated integrals 
for double loop spaces.  hep-th/9403055.

\bibitem{Hirs03} P. Hirschhorn, {\em Model categories and their localizations.} Mathematical 
Surveys and Monographs, {\bf 99}. American Mathematical Society, Providence, RI,  2003.

 \bibitem{Lill73} J. Lillig,  A union theorem for cofibrations.  
 {\em Arch. Math. (Basel)}  {\bf 24}  (1973), 410--415.
 

\bibitem{Kont99} M. Kontsevich, Operads and motives in deformation
quantization, {\em{Lett. Math. Phys.}} 48 (1999), 35--72.

\bibitem{Kriz95} I. Kriz, On the rational homotopy type of configuration spaces.  
{\em Ann. of Math. (2)} {\bf 139}  (1994),  no. 2, 227--237. 

\bibitem{KuTh01} G. Kuperberg and D. Thurston,   Perturbative 3-manifold invariants by 
cut-and-paste topology.  math.GT/9912167.

\bibitem{Mark99} M. Markl,   A compactification of the real configuration space as an 
operadic  completion.  {\em  J. Algebra}  {\bf 215}  (1999),  no. 1, 185--204.

\bibitem{MSS02} M. Markl, S. Shnider, and J. Stasheff, 
{\em Operads in algebra, topology and physics}. Math. Surv. and Monographs, 96. 
AMS, Providence, RI, 2002.

\bibitem{Poir02} S. Poirier, The configuration space integral for links in $\R\sp 3$.  
{\em Algebr. Geom. Topol.} {\bf 2} (2002), 1001--1050.

\bibitem{Sane01}  S. Saneblidze and R. Umble,  A diagonal on the
associahedra. math.AT/0011065 (2001). 

\bibitem{CMS03} D. Sinha, Algebraic and differental topology of configuration spaces, in preparation.

\bibitem{McSS03} D. Sinha, Operads and knot spaces.  math.AT/0407039.

\bibitem{Sinh02}  D. Sinha,  The topology of space of knots.  math.AT/0202287.

\bibitem{Stas63} J. Stasheff, Homotopy associativity of $H$-spaces, I. {\em 
Trans. Amer. Math. Soc.} {\bf 108} (1963) 275--292.

\bibitem{Tota96} B. Totaro, Configuration spaces of
algebraic varieties. {\em Topology} {\bf 35} (1996), no. 4,
1057--1067.

\bibitem{Yuzv97} S. Yuzvinsky, Cohomology bases for the De Concini-Procesi models of hyperplane arrangements and sums over trees.  {\em Invent. math.} {\bf 127}, 
(1997), 319--335.

\end{thebibliography}
\end{document}